\documentclass[times]{article}

\usepackage{fullpage}
\usepackage{moreverb}
\usepackage{tikz}
\usepackage{amsmath,amssymb,amsbsy}
\usepackage{bm,multirow}
\usepackage[dvips,colorlinks,bookmarksopen,bookmarksnumbered,citecolor=black,urlcolor=black,breaklinks=true]{hyperref}
\usepackage{booktabs}
\usepackage{mathrsfs}
\usepackage{algorithm,algorithmic}
\usepackage{setspace}

\usetikzlibrary{positioning,calc, shapes, backgrounds}

\graphicspath{{../ps/}}
\tikzset{>=latex}

\def\equationautorefname~#1\null{%
  Eq.~(#1)\null
}
\def\algorithmautorefname~#1\null{%
  Algorithm~#1\null
}
\def\sectionautorefname~#1\null{%
  Section~#1\null
}
\def\subsectionautorefname~#1\null{%
  Subsection~#1\null
}

\newcommand\mO{\mathcal{O}}
\newcommand\mH{\mathcal{H}}

\newcommand{\bd}{\partial}
\newcommand{\ff}{\mathbf{f}}

\newcommand{\grad}{\triangledown}
\newcommand{\nn}{\mathbf{n}}
\newcommand{\rr}{\mathbf{r}}

\newcommand{\ssigma}{\boldsymbol{\sigma}}
\newcommand{\uu}{\mathbf{u}}
\newcommand{\xx}{\mathbf{x}}
\newcommand{\yy}{\mathbf{y}}

\title{An efficient preconditioner for the fast simulation of a 2D Stokes flow in porous media}
\author{Pieter Coulier, Bryan Quaife, and Eric Darve}

\begin{document}
\maketitle

\begin{abstract}
We consider an efficient preconditioner for boundary integral equation
(BIE) formulations of the two-dimensional Stokes equations in porous
media.  While BIEs are well-suited for resolving the complex porous
geometry, they lead to a dense linear system of equations that is
computationally expensive to solve for large problems.  This expense is
further amplified when a significant number of iterations is required in
an iterative Krylov solver such as GMRES.  In this paper, we apply a
fast inexact direct solver, the inverse fast multipole method (IFMM), as
an efficient preconditioner for GMRES. This solver is based on the
framework of $\mH^{2}$-matrices and uses low-rank compressions to
approximate certain matrix blocks. It has a tunable accuracy
$\varepsilon$ and a computational cost that scales as $\mO (N \log^2
1/\varepsilon)$. We discuss various numerical benchmarks that validate
the accuracy and confirm the efficiency of the proposed method. We
demonstrate with several types of boundary conditions that the
preconditioner is capable of significantly accelerating the convergence
of GMRES when compared to a simple block-diagonal preconditioner,
especially for pipe flow problems involving many pores.
\end{abstract}

\section{Introduction}\label{sec:introduction}
\begin{figure}[hbtp!]
  \begin{center}
  \includegraphics*[width=0.975\linewidth,clip=true,bb=0 8 238 39]{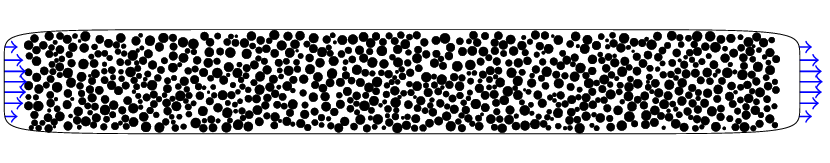}
  \end{center}
  \caption{Geometry of the problem: circular pores surrounded by an
  outer boundary.  The quiver plot on the intake and outtake shows the
  pipe flow boundary conditions defined in \autoref{eq:noslip_bcs}.}
  \label{fig:schematic}
\end{figure}

Flow inside a porous medium finds many applications in natural and
engineering systems.  Subsurface flows~\cite{cue-jua2008,
kit-jia-val-chi-tsu-chr2014}, erosion~\cite{mcd-hun-sit1986,
var-sta-pap1996}, fuel cells~\cite{wan2004}, and filtration
systems~\cite{lec-bot-wie2004, wu-liu-xie-liu-sun2012} are a few
examples of physical processes that are governed by low Reynolds number
porous media flow.  Accurately resolving such flows is essential for
modelling transport~\cite{wu-ye-sud2010}, mixing~\cite{leb-den-vil2013,
hid-fe-cue-jua2012}, anomalous diffusion~\cite{ber-sch-sil2000,
sun-che-che2009}, breakthrough curves~\cite{hag-mck-mei2000}, carbon
dioxide sequestration~\cite{kaz-blo-kyr-chr2015, szu-hes-jua2013},
uncertainty quantification~\cite{ica-boc-tem2016, mey-tch-jen2013},
chemical reactions~\cite{den-leb-eng-bij2011, tel-dem-gre-kar2005}, and
many other applications.  In most of these works, in particular, those
that focus on numerics and modelling, the porous medium is assumed to be
two dimensional.  Therefore, there is a strong need for efficient and
accurate numerical methods to simulate complex two-dimensional porous
media flow.  In this paper, we focus on solving the two-dimensional
incompressible Stokes equations in the geometry illustrated in
\autoref{fig:schematic}.

The incompressible Stokes equations have been solved in porous domains
resembling \autoref{fig:schematic} using finite
differences~\cite{leb-oge-leb-hul-mar-sch1996}, finite
volumes~\cite{ica-boc-mar-tos-set2014,and-alm-fil-hav-suk-sta1997},
algebraic multigrid~\cite{blu-bij-don-gha-igl-mos-pal-pen2013},
Lattice-Boltzmann~\cite{man-gla-war1999}, and continuous time random
walks~\cite{dea-leb-den-tar-bol-dav2013}.  However, each of these
methods struggle with either enforcing incompressibility, meshing the
geometry, or obtaining high-order accuracy.  Boundary integral equations
(BIEs), on the other hand, are a powerful method that addresses each of
these challenges.  In particular, BIEs are naturally adaptive, allow for
high-order approximations, automatically enforce the incompressibility
constraint, and enable the computation of the velocity field without
computing the pressure.  Another advantage of a BIE is that the unknown
two-dimensional velocity field is written in terms of a one-dimensional
integral over the boundary of the domain.  This dimension reduction
implies that only the boundary, which is a collection of one-dimensional
closed curves, needs to be meshed.  Furthermore, a natural extension of
this work is to solve the three-dimensional Stokes equations while only
discretizing the two-dimensional boundary of the geometry.  Other groups
have used integral equations to simulate porous media by using a Darcy
approximation~\cite{liu-che-lig1981, lou-lee-kam1998}, or by coupling a
penalization method with a volume integral
equation~\cite{mal-gho-bir2014}.  In contrast, we will use a boundary
integral equation to solve the incompressible Stokes equations.

The main disadvantage of BIEs is the need to solve a dense linear system
of equations.  A direct solver based on Gaussian elimination requires
$\mO(N^{3})$ operations---where $N$ denotes the number of
unknowns---which is impractical for realistic geometries. It also
requires the assembly and storage of the dense matrix, leading to an
$\mO(N^{2})$ memory requirement. Iterative solvers such as GMRES, on the
other hand, only rely on (accelerated) matrix-vector multiplications and
have an $\mO(N)$ memory requirement, but can suffer from a poor
convergence behavior due to the ill-conditioning of the linear system.
In the recent decade, there have been two main directions of research
that have lead to faster algorithms for solving discretized BIEs.

A first approach focuses on developing integral equations that result
in linear systems that are both well-conditioned and require a
mesh-independent number of GMRES iterations to reach
convergence~\cite{cam-ips-kel-mey-xue1996}.  These properties generally
follow from second-kind integral equation methods with compact integral
operators.  Because of the mesh-independent number of GMRES iterations,
the asymptotic CPU time scales with the cost of a single dense
matrix-vector multiplication.  Using fast summation methods such as the
fast multipole method (FMM)~\cite{gre-rok1987, yin-bir-zor2004}, Ewald
summation~\cite{dar-yor-ped1993}, or tree-code
methods~\cite{bar-hut1986}, the cost of a matrix-vector multiplication
can be reduced from $\mO(N^{2})$ to either $\mO(N\log N)$ or $\mO(N)$.
While this approach can result in an optimal asymptotic complexity,
complexities of the geometry, such as nearly touching pores and regions
of high curvature, often still results in a large number of GMRES
iterations~\cite{qua-bir2015a}.  In addition, the double-layer
potential formulation, which requires a mesh-independent number of
GMRES iterations, has less regularity than the single-layer
potential~\cite{col-kre2012, hsi-wen2008}, which requires a
mesh-dependent number of GMRES iterations. To reduce the number of
iterations, a preconditioner can be incorporated.  While simple
techniques such as block-diagonal preconditioners are easy to
construct, they are not always capable of reducing the number of
iterations significantly.  An incomplete list of more sophisticated
preconditioners, which vary in construction time and efficacy, include
additive or multiplicative Schwarz~\cite{heu:ste:tra1998,
mai-ste-tra2000}, BPX~\cite{fun-ste1997}, sparse approximate
inverses~\cite{che2000, gra-kum-sam1998}, approximate LU
decompositions~\cite{beb2005}, clustering~\cite{nab-kor-lei-whi1994,
vav1992}, integral equations of opposite order~\cite{ste-wen1998},
Calderon identities~\cite{hsi-wen2008}, or
multigrid~\cite{bra-ley-pas1994, hsi-xu-zha2014, lan-pus-rei2003,
of2008}.  In the case of multigrid preconditioners, specialized
smoothers must be developed since standard smoothers like Jacobi or
Gauss-Seidel will smooth the low, rather than the high, frequency
components.  Examples of such smoothers may require developing
approximate inverses of the integral operator~\cite{bra-ley-pas1994}
involving, for example, the Laplace-Beltrami
operator~\cite{lan-pus2005}.

A second approach targets the use of fast direct solvers that exploit
the structure of the linear system that arises when discretizing BIEs.
The structure can be loosely described as matrices with hierarchical
off-diagonal blocks being numerically low-rank.  These matrices are
formalized as, with an increasing level of complexity, hierarchically
off-diagonal low-rank (HODLR) matrices~\cite{amin14a}, hierarchically
semi-separable (HSS) matrices~\cite{chan06a,shen07b},
$\mH$-matrices~\cite{beb2005,bebe08a}, and
$\mH^{2}$-matrices~\cite{borm2007data,hackbusch2002data}. Fast methods
for efficiently constructing and storing approximate inverses of these
matrices have been developed in recent years~\cite{gill14a, kong11a,
mart05b}, mainly for the aforementioned HODLR and HSS matrices.
Moreover, once the approximate inverse, which only depends on the
geometry, is computed and stored, it can be applied with a minimal
amount of CPU time to multiple right-hand sides.  Therefore, direct
solvers are extremely useful when applying a temporal discretization to
a partial differential equation~\cite{mar-bar-gil-vee2015}, or when
considering scattering for many different incident
waves~\cite{gil-bar2013, gil-bar-mar2014}.  


We blend the two aforementioned approaches in this paper: we present the
inverse fast multipole method (IFMM)~\cite{amb-dar2014} as an inexact
fast direct solver based on $\mH^{2}$-matrices and apply it as a
preconditioner in a fast multipole accelerated GMRES
solver~\cite{ij-sjsc-coul-a}. The preconditioner has a tunable accuracy
$\varepsilon$ and its computational cost scales almost linearly with the
problem size as $\mO (N \log^2 1/\varepsilon)$.  We will demonstrate
that the IFMM is capable of significantly reducing the number of
iterations and the CPU time.

The outline of this paper is as follows. \autoref{sec:StokesFlow}
discusses the fluid model and layer potential formulation.  The inverse
fast multipole method (IFMM) is subsequently introduced in
\autoref{sec:IFMM}, focusing on the main features for achieving an
efficient preconditioner.  Numerical benchmarks are presented in
\autoref{sec:NumericalExamples}.  We report results for both the
unpreconditioned system, and for preconditioned systems where
block-diagonal and IFMM preconditioners (at several accuracies) are
used. Summarizing conclusions are finally drawn in
\autoref{sec:Conclusions}.

\section{Stokes flow in porous media}\label{sec:StokesFlow}
We are interested in simulating an incompressible Newtonian fluid in the
geometry illustrated in \autoref{fig:schematic}.  We consider scales
where the dimensionless Reynolds number is small so that the fluid is
assumed to be Stokesian.  If we assumed a homogeneous porosity, the much
simpler Darcy equations could be solved.  However, we do not make this
assumption and we focus on accurately and efficiently solving the
incompressible Stokes equation with a no-slip boundary condition on each
pore.  We also assume that the flow is two-dimensional and that there
are no external body forces such as gravity.

We avoid the Stokes paradox by bounding the porous region by the boundary
$\Gamma_{0}$.  To resemble a pipe flow, we let $\Gamma_{0}$ be a
mollification of the rectangle $[0,L] \times [-H,H]$ (see
\autoref{fig:schematic}).  By smoothing the corners of $\Gamma_{0}$,
specialized quadrature for integral equations is avoided.  Each inner
boundary, denoted by $\Gamma_{k}$, $k=1,\ldots,M$, is a circle of
variable radii.  The area enclosed by the outer boundary is denoted by
$\Omega_{0}$, and the area enclosed by each pore by $\Omega_{k}$.
Therefore, the geometry is given by
\begin{align*}
  \Omega = \Omega_{0} \setminus \bigcup_{k=1}^{M} \Omega_{k},
\end{align*}
and its boundary is $\Gamma = \Gamma_{0} \cup \Gamma_{1} \cup \cdots
\cup \Gamma_{M}$.  

With this setup, by defining the fluid velocity $\uu(\xx)$, the pressure
$p(\xx)$, the viscosity $\mu$, and a prescribed flow on the boundary
$\ff(\xx)$, the governing equations are 
\begin{equation}
\begin{aligned}
  \mu \Delta \uu(\xx) &= \grad p(\xx), \: \grad \cdot \uu(\xx) = 0, \quad 
    &&\xx \in \Omega, \\
  \uu(\xx) &= \ff(\xx), &&\xx \in \bd \Omega.
  \label{eq:stokes}
\end{aligned}
\end{equation}
We impose boundary conditions that correspond to a pipe flow at the
intake and outtake, and no slip on the boundary of the pores
\begin{align}
  \ff(\xx) = \left\{
  \begin{array}{ll}
      k(H^{2} - y^{2}), & \xx \in \Gamma_{0}, \\
      0, & \xx \in \displaystyle\bigcup_{k=1}^{M} \Gamma_{k},
  \end{array}
  \right.
  \label{eq:noslip_bcs}
\end{align}
where $k$ sets the velocity scale.  Furthermore, we assume that $\mu=1$
from this point onwards.

\subsection{Integral equation formulation}
A BIE formulation of \autoref{eq:stokes} has several advantages over the
differential form.  The pressure need not be computed, the
incompressibility constraint is automatically satisfied, and the
velocity is written in terms of a layer potential involving an unknown
density function defined only on $\Gamma$.  Two popular choices for the
layer potential are the single-layer and double-layer potentials
\begin{alignat*}{3}
  \mathcal{S}[\ssigma](\xx) &= \frac{1}{4\pi} \int_{\Gamma} \left( 
  -\log \rho \mathbf{I} + \frac{\rr \otimes \rr}{\rho^{2}}\right) 
  \ssigma(\yy) ds_{\yy}, \qquad && \xx \in \Omega, \\
  \mathcal{D}[\ssigma](\xx) &= \frac{1}{\pi}\int_{\Gamma}
  \frac{\rr \cdot \nn}{\rho^{2}} \frac{\rr \otimes \rr}{\rho^{2}}
  \ssigma(\yy) ds_{\yy}, && \xx \in \Omega,
\end{alignat*}
respectively, where $\rr = \xx - \yy$, $\rho = \|\rr\|$, $\nn$ is the
unit outward normal, and $\ssigma$ is a density function defined on
$\Gamma$.  Given our Dirichlet boundary condition, the double-layer
potential formulation results in a second-kind integral equation whose
condition number is guaranteed to be mesh-independent.  However, the
single-layer potential requires less smoothness of $\ssigma$
since~\cite{col-kre2012, hsi-wen2008}
\begin{align*}
  \mathcal{S}: H^{-\frac12}(\Gamma) &\rightarrow H^{1}(\Omega), \\
  \mathcal{D}: H^{\frac12}(\Gamma) &\rightarrow H^{1}(\Omega),
\end{align*}
are continuous.  Moreover, this result holds even under the relaxed
assumption that $\Omega$ is a bounded Lipschitz domain.  Another
advantage of the single-layer potential is that it has full-rank, as
opposed to the double-layer potential which has a rank three null
space~\cite{pow1993} for each connected component of $\Gamma$.  While
this rank deficiency can be removed, it is unclear if this is
compatible with the IFMM.  Moreover, for complex geometries, the
double-layer potential may still require many GMRES
iterations~\cite{qua-bir2015a}, but in this case, we expect that
applying the IFMM preconditioner will result in a computational speedup.

Given the additional regularity and full-rank of the single-layer
potential, we choose to represent the velocity $\uu$ as
\begin{align}
  \uu(\xx) = \mathcal{S}[\ssigma](\xx) = \frac{1}{4\pi} \int_{\Gamma} \left( 
  -\log \rho \mathbf{I} + \frac{\rr \otimes \rr}{\rho^{2}}\right) 
  \ssigma(\yy) ds_{\yy}, \quad \xx \in \Omega,
  \label{eq:stokesSLP}
\end{align}
where $\ssigma$ is an unknown density function.  Taking the limiting
value of \autoref{eq:stokesSLP} as $\xx$ tends to $\Gamma$, the density
function $\ssigma$ must satisfy the first-kind integral
equation~\cite{poz1992}
\begin{align}
  \ff(\xx) = \frac{1}{4\pi} \int_{\Gamma} \left( -\log \rho \mathbf{I} +
  \frac{\rr \otimes \rr}{\rho^{2}}\right) \ssigma(\yy) ds_{\yy}, 
  \quad \xx \in \Gamma.
  \label{eq:stokesBIE}
\end{align}
It is known that upon discretizing \autoref{eq:stokesBIE}, the result is
an ill-conditioned linear system that requires a mesh-dependent number
of GMRES iterations.  However, the linear system is amenable to the
solvers for structured matrices.  Once the unique solution of
\autoref{eq:stokesBIE} is computed, the density function $\ssigma$ is
substituted into \autoref{eq:stokesSLP} to evaluate the velocity
$\uu(\xx)$ at any point $\xx \in \Omega$.

\subsection{Discretization of the BIE}
We adopt a collocation method to discretize \autoref{eq:stokesBIE}.
First, each interior curve $\Gamma_{k}$, $k=1,\ldots,M$ is parameterized
and discretized at $N_{\mathrm{int}}$ points, and the outer wall is
discretized at $N_{\mathrm{ext}}$ points.  The resulting discretization
points are $\xx_{i}$, $i=1,\ldots,MN_{\mathrm{int}} + N_{\mathrm{ext}}$.
We enforce the BIE at these $N = MN_{\mathrm{int}} + N_{\mathrm{ext}}$
collocation points by requiring
\begin{align*}
  \ff(\xx_{i}) = \frac{1}{4\pi} \int_{\Gamma} \left( -\log
  \|\xx_{i}-\yy\| \mathbf{I} +
  \frac{(\xx_{i} - \yy) \otimes (\xx_{i} - \yy)}{\|\xx_{i} - \yy\|^{2}}\right)  
  \ssigma(\yy) ds_{\yy}, \quad i=1,\ldots, N.
\end{align*}
Then, quadrature is applied, where the quadrature nodes coincide with
the discretization points $\xx_{i}$.  This results in the dense linear
system
\begin{align}
  \ff(\xx_{i}) = \sum_{j=1}^{N} w_{j}G(\xx_{i},\xx_{j})
    \ssigma(\xx_{j}), \quad i=1,\ldots N, \label{eq:dense}
\end{align}
where $G(\xx,\yy)$ is the kernel of the integral operator in
\autoref{eq:stokesSLP}.  The weights $w_{j}$ are the product of the
arclength term with appropriate quadrature weights.  

We integrate the weak logarithmic singularity by applying a sixth-order
quadrature rule of Kapur and Rokhlin~\cite{kapu97a}.  This quadrature
rule is identical to the trapezoid rule, except that the weights are
modified at the six nodes to the left and right of the singularity.
Therefore, the quadrature nodes coincide with the discretization points
$\xx_{i}$, $i=1,\ldots,N$.  Since there are only $\mO(1)$ quadrature
weights that are modified for each $\xx_{i}$, these quadrature rules are
compatible with the FMM.  In addition, in contrast to Alpert quadrature
rules~\cite{alp1999}, additional quadrature nodes are not introduced, so
no interpolation is required.

Given this quadrature formula,~\autoref{eq:dense} results in a dense
$2N \times 2N$ linear system that, for brevity, we write as
\begin{align}
  \ff = \mathbf{A} \ssigma.
  \label{eq:dense_matrix}
\end{align}
The matrix $\mathbf{A}$ is a structured matrix with numerical low-rank
off-diagonal blocks.  GMRES can be used to iteratively
solve~\autoref{eq:dense_matrix}, and each matrix-vector multiplication
can be done with $\mO(N)$ complexity by using the FMM.  However, the
number of required GMRES iterations will be large, and depend on the
mesh size.  Therefore, we will use the structure of $\mathbf{A}$ to
develop an efficient preconditioner for GMRES.

Upon computing the density function $\ssigma$ at the points
$\xx_{1},\ldots,\xx_{N}$, the velocity field needs to be computed at
points $\xx \in \Omega$.  This can be done with the trapezoid rule
\begin{align}
  \uu(\xx) \approx \frac{1}{4\pi} \sum_{j=1}^{N} 
    G(\xx,\xx_{j})\ssigma(\xx_{j}) \Delta s_{j},  
    \label{eq:stokesSLP_discretized}
\end{align}
where $\Delta s_{j}$ is an arclength spacing of $\Gamma$ at $\xx_{j}$.
The trapezoid rule guarantees spectral accuracy, but this asymptotic
convergence rate is delayed when $\xx$ is too close to $\Gamma$\footnote{Assuming the geometry is resolved and the density
function is exact, the trapezoid rule results in machine precision if
$d(\xx,\Gamma) > \sqrt{\Delta s}$.}.  There are methods to guarantee a
uniform error for all $\xx \in
\Omega$~\cite{bar-wu-vee2014,yin-bir-zor2006}, but this is not the
focus of this work.  Instead, we focus on efficiently solving
\autoref{eq:dense_matrix} for the density function $\ssigma$.

\section{The inverse fast multipole method (IFMM)}
\label{sec:IFMM}
The inverse fast multipole method (IFMM)~\cite{amb-dar2014} is an
inexact\footnote{The error can be controlled and made as small as
needed, however.} fast direct solver for $\mH^{2}$-matrices that can be
applied as a preconditioner for GMRES.  It has successfully been used to
precondition matrices arising from integral
operators~\cite{ij-sjsc-coul-a} and radial basis
interpolation~\cite{cou-dar2016}.  The main features of the algorithm
are concisely described in the following subsections; the reader is
referred to \cite{ij-sjsc-coul-a} for a more detailed description.
Although a 2D problem is considered in this paper, the IFMM is
applicable to 3D problems as well (see~\cite{ij-sjsc-coul-a} for
examples).

\subsection{$\mH^{2}$-matrices}
\label{subsec:H2}
A variety of fast direct solvers for dense linear systems such as
\autoref{eq:dense_matrix} have been developed in recent years.  These
methods are based on a multilevel, hierarchical decomposition of the
physical domain of interest, from the root down to the leaf level~$l$
(e.g., by means of a quadtree  in 2D or an octree in 3D, either uniform
or adaptive).  This is subsequently exploited to represent the matrix
$\mathbf{A}$ in a compressed format. 

The simplest of these methods assume that all off-diagonal blocks are
low-rank (only self-interactions are considered to be of full-rank),
leading to so-called hierarchically off-diagonal low-rank
(HODLR)~\cite{amin14a} and hierarchically semi-separable (HSS)
matrices~\cite{chan06a,shen07b}. In the latter, a nested approach is
used (i.e.,~the low-rank basis at a certain hierarchical level is
constructed using the low-rank basis at the child level), while this is
not the case in the former. Given this assumption on the matrix
structure, the inversion can be performed exactly. Direct $\mO(N)$
algorithms for HODLR matrices can be found in~\cite{amin14a,kong11a},
while fast solvers for HSS matrices have been presented by, among
others, Martinsson and Rokhlin~\cite{mart05b}, Chandrasekaran et
al.~\cite{chan06a}, Xia et al.~\cite{xia10a,xia10b}, and Gillman et
al.~\cite{gill14a}.

The assumption that all off-diagonal blocks are low-rank is a major
drawback for the aforementioned HODLR and HSS matrices. Indeed, the rank
is in that case likely to grow in an unbounded fashion for 2D and 3D
problems with a complex geometry~\cite{ambi13b}, hence jeopardizing the
computational efficiency of the algorithms. This drawback can be
circumvented by adopting the more general framework of $\mH$- or
$\mH^{2}$-matrices (non-nested vs.~nested), in which only
non-neighboring interactions are replaced by low-rank approximations.
This is the approach followed in the IFMM~\cite{amb-dar2014}. It is
worth noting that the algorithm presented in~\cite{mart05b} has recently
been extended by Corona et al.~\cite{coro15a} to obtain a direct solver
for HSS matrices that scales for 2D problems as $\mO(N)$.  We note that
since our problem is two-dimensional, a HSS approximation of
\autoref{eq:dense_matrix} would potentially have bounded rank.  However,
to extend the work to three dimensions, we choose to use an $\mH^{2}$
framework.

An $\mH^{2}$-matrix $\mathbf{A}^{(l)}_{\mH}$ is weak hierarchical matrix
(in which only non-neighboring, well-separated interactions are assumed
to be low-rank) with a nested basis that approximates the original
matrix $\mathbf{A}$. It can be expressed as~\cite{borm2007data,
hackbusch2002data}
\begin{equation}
    \mathbf{A} \simeq 
    \mathbf{A}^{(l)}_{\mH} = \mathbf{S}^{(l)} + 
      \mathbf{U}^{(l)} \mathbf{A}^{(l-1)}_{\mH}
      \mathbf{V}^{(l)^\mathrm{T}}.
    \label{eq:A_H}
\end{equation}
$\mathbf{A}^{(l)}_{\mH}$ is thus decomposed into a block sparse matrix,
$\mathbf{S}^{(l)}$, containing all the neighboring and self-interactions
at the leaf level $l$, and a low-rank term, $\mathbf{U}^{(l)}
\mathbf{A}^{(l-1)}_{\mH} \mathbf{V}^{(l)^\mathrm{T}}$, that
characterizes all well-separated interactions. $\mathbf{U}^{(l)}$ and
$\mathbf{V}^{(l)}$ are interpolation and anterpolation operators,
respectively, while $\mathbf{A}^{(l-1)}_{\mH}$ is an $\mH^{2}$-matrix of
a smaller size. \autoref{eq:A_H} is applied recursively to
$\mathbf{A}^{(l-1)}_{\mH}$ until the top level\footnote{The recursion
stops at $\mathbf{A}^{(1)}_{\mH}$.}, where there are no more
well-separated interactions, is reached. The IFMM is a novel algorithm
for the approximate inversion of an $\mH^{2}$-matrix, and consists of
two crucial steps: extended sparsification and compression of fill-ins
that appear throughout the elimination.

\subsection{Extended sparsification}
\label{subsec:sparsification}
The hierarchical structure of $\mathbf{A}^{(l)}_{\mH}$ can be exploited
to represent $\mathbf{A}^{(l)}_{\mH} \ssigma = \ff$ as an extended
sparse system, which is more attractive from a computational point of
view. Introducing auxiliary variables $\mathbf{y}^{(l)} =
\mathbf{V}^{(l)^\mathrm{T}} \ssigma$ and $\mathbf{z}^{(l)} =
\mathbf{A}^{(l-1)}_{\mH} \mathbf{y}^{(l)}$ (corresponding to multipole
and local coefficients in the FMM) leads to the following extended
system:
\begin{equation}
	\begin{bmatrix}
		\mathbf{S}^{(l)} 				& \mathbf{U}^{(l)}			& \\
		\mathbf{V}^{(l)^\mathrm{T}} 	& 			             	& -\mathbf{I}\\ 
										& -\mathbf{I}				& \mathbf{A}^{(l-1)}_{\mH}
	\end{bmatrix}
	\begin{bmatrix}
		\ssigma\\
		\mathbf{z}^{(l)}\\
		\mathbf{y}^{(l)}
	\end{bmatrix}
	=
	\begin{bmatrix}
		\ff\\
		\mathbf{0}\\
		\mathbf{0}
	\end{bmatrix}.
	\label{eq:extendedsparse1}
\end{equation}
In \autoref{eq:extendedsparse1}, the sparse blocks $\mathbf{S}^{(l)}$,
$\mathbf{U}^{(l)}$, and $\mathbf{V}^{(l)^\mathrm{T}}$ are situated in
the top left corner of the matrix, while the remaining dense block
$\mathbf{A}^{(l-1)}_{\mH}$ is pushed to the bottom right corner of the
matrix. As the latter is an $\mH^{2}$-matrix as well, the same
technique can be applied recursively to further extend and sparsify the
system.  One finally obtains the following system of equations, which
is equivalent to $\mathbf{A}^{(l)}_{\mH} \ssigma = \ff$:
\begin{equation}
	\begin{bmatrix}
		\mathbf{S}^{(l)} 				& \mathbf{U}^{(l)}			& & & & & & & \\
		\mathbf{V}^{(l)^\mathrm{T}} 	&             				& -\mathbf{I}& & & & & & \\ 
				 						& -\mathbf{I}				& \mathbf{S}^{(l-1)} & \mathbf{U}^{(l-1)} & & & & & \\
				 						& 							& \mathbf{V}^{(l-1)^\mathrm{T}}	 &  & -\mathbf{I}& & & & \\
				 						& 							& 	 & -\mathbf{I} & \mathbf{S}^{(l-2)} & & & & \\
				 						& 							& 	 &  & & \ddots & & & \\
				 						& 							& 	 &  & & & \mathbf{S}^{(2)} & \mathbf{U}^{(2)} & \\
				 						& 							& 	 &  & & & \mathbf{V}^{(2)^\mathrm{T}} & & -\mathbf{I}\\
				 						& 							& 	 &  & & & &-\mathbf{I} & \mathbf{A}^{(1)}_{\mH}
	\end{bmatrix}
	\begin{bmatrix}
		\ssigma\\
		\mathbf{z}^{(l)}\\
		\mathbf{y}^{(l)}\\
		\mathbf{z}^{(l-1)}\\
		\mathbf{y}^{(l-1)}\\
		\vdots\\
		\mathbf{y}^{(3)}\\
		\mathbf{z}^{(2)}\\
		\mathbf{y}^{(2)}
	\end{bmatrix}
	=
	\begin{bmatrix}
		\ff\\
		\mathbf{0}\\
		\mathbf{0}\\
		\mathbf{0}\\
		\mathbf{0}\\
		\vdots\\
		\mathbf{0}\\
		\mathbf{0}\\
		\mathbf{0}
	\end{bmatrix}.
	\label{eq:extendedsparse2}
\end{equation}
The variables and equations are ordered in such a way that a sparse
matrix with a symmetric fill-in pattern is obtained. Note that this
extended system is only moderately larger than the original system,
since the dimensions of the auxiliary variables (characterizing the
low-rank interactions) are small compared to the dimension of
$\ssigma$.

\subsection{Compression of fill-ins}
\label{subsec:compression}
The application of Gaussian elimination to the extended sparse system of
\autoref{eq:extendedsparse2} leads to the creation of dense fill-in
blocks, which jeopardizes the computational efficiency of the method. In
order to obtain a fast solver, it is crucial to preserve the sparsity
pattern throughout the elimination. This can be achieved by compressing
and redirecting fill-ins that correspond to well-separated interactions,
as these are expected to have a numerically low-rank. 

For example, consider the fill-ins $\mathbf{S}^{(l)\prime}_{ij}$ and
$\mathbf{S}^{(l)\prime}_{ji}$ arising between the well-separated
variables $\ssigma_i$ and $\ssigma_j$ (the prime indicates a fill-in;
$\mathbf{S}^{(l)}_{ij}$ and $\mathbf{S}^{(l)}_{ji}$ are originally zero
blocks). These fill-ins can be approximated by a low-rank representation
(compression step) by means of a truncated singular value decomposition
(SVD) in which only the most significant singular values and vectors are
retained:
\begin{align}
  \mathbf{S}^{(l)\prime}_{ij} &\simeq \mathbf{U}^{(l)\prime}_{i}
  \mathbf{\Sigma}^{(l)\prime}_{ij}
  \mathbf{V}^{(l)\prime^\mathrm{T}}_{j}\label{eq:S_ij},\\
  \mathbf{S}^{(l)\prime}_{ji} &\simeq \mathbf{U}^{(l)\prime}_{j}
  \mathbf{\Sigma}^{(l)\prime}_{ji}
  \mathbf{V}^{(l)\prime^\mathrm{T}}_{i}\label{eq:S_ji}.
\end{align}
The rank of approximations \autoref{eq:S_ij} and \autoref{eq:S_ji} is
determined by a prescribed tolerance~$\varepsilon$ (either relative or
absolute) on the singular values.

Next, a recompression step is performed to obtain new low-rank
interpolation and anterpolation operators
$\widehat{\mathbf{U}}^{(l)}_{i}$ and $\widehat{\mathbf{V}}^{(l)}_{i}$,
respectively, as well as $\widehat{\mathbf{U}}^{(l)}_{j}$ and
$\widehat{\mathbf{V}}^{(l)}_{j}$, such that the fill-ins
$\mathbf{S}^{(l)\prime}_{ij}$ and $\mathbf{S}^{(l)\prime}_{ji}$ can be
redirected through the existing low-rank interaction between $\ssigma_i$
and $\ssigma_j$ in $\mathbf{A}^{(l-1)}_{\mH}$. As a result, there is no
need to store $\mathbf{S}^{(l)\prime}_{ij}$ and
$\mathbf{S}^{(l)\prime}_{ji}$ explicitly, implying that the sparsity
pattern of the matrix is maintained throughout the elimination. If we
make the assumption that the rank remains bounded as the problem size
increases, it can be demonstrated that this leads to an algorithm with
an asymptotic complexity of $\mO (N \log^2 1/\varepsilon)$. Although
there is no mathematical proof to support this assumption, numerical
examples confirm the fact that this is indeed the case in many practical
examples~\cite{ij-sjsc-coul-a,lize14a}.

The compression and recompression procedures make the IFMM inexact, but
the accuracy can be tuned by varying $\varepsilon$. In terms of
computational efficiency, one needs to make a trade-off  between a
highly accurate direct solver or a low-accuracy preconditioner. We focus
in this paper on its use as a preconditioner. More details on the
aforementioned procedures are provided in~\cite{ij-sjsc-coul-a}. The
value of $\varepsilon$ should be chosen in relation to the accuracy of
the initial low-rank operators in \autoref{eq:A_H} (i.e., the higher the
initial rank, the smaller $\varepsilon$ should be).

\section{Numerical benchmarks}
\label{sec:NumericalExamples}

Several numerical benchmarks are considered in this section to
demonstrate the efficiency of the IFMM as a preconditioner for the
problems under concern. In all of these benchmarks,
\autoref{eq:dense_matrix} is solved for the density vector
$\mathbf{\boldsymbol{\sigma}}$ with a fast multipole accelerated GMRES
solver~\cite{saa-sch1986}, and \autoref{eq:stokesSLP_discretized} is
subsequently used to evaluate the velocities $\mathbf{u}(\mathbf{x})$
inside the domain. A tolerance of $10^{-8}$ is specified for the
relative residual in GMRES.  The computations have been performed on
Intel\textsuperscript{\textregistered}
Xeon\textsuperscript{\textregistered} E5-2650 v2 (2.60~GHz) CPUs.

In all benchmarks, a black-box approach is used in the IFMM for
constructing the initial low-rank operators in \autoref{eq:A_H}. More
precisely, interpolation based on Chebyshev polynomials \cite{fong09a}
is employed (with $n$ Chebyshev nodes in each direction), followed by an
additional SVD to reduce the rank. A uniform quadtree decomposition of
the domain is used, and the number of levels is adjusted to provide a
reasonable trade-off between the time spent at the leaf level and the
time spent at higher levels in the tree\footnote{Increasing the number
of levels makes the elimination of the leaf nodes faster, but increases
the time spent in the tree, and vice versa.}.

We are interested in the porous geometry illustrated in
\autoref{fig:schematic}.  The outer wall is a rounded off version of the
rectangle $[0,42] \times [-2.6,2.6]$, and the radii of the pores are
distributed in $[6.4 \times 10^{-2},2.5 \times 10^{-1}]$.  The smallest
gap between two pores is $4.7 \times 10^{-3}$, and the smallest gap
between a pore and the outer wall is $2.7 \times 10^{-2}$.  Finally, the
porous region has a porosity of $55.4\%$.

We test the preconditioner on simpler geometries by considering only the
first $22$ and $226$ left-most pores.  We also scale the outer boundary
so that the area to the right of the pores remains roughly constant.
Before studying the pipe flow boundary conditions defined in
\autoref{eq:noslip_bcs}, we first validate our method by examining the
convergence with two sets of boundary conditions that have a closed-form
solution in \autoref{subsec:shear} and~\ref{subsec:analytical}.  We
subsequently discuss the pipe flow boundary conditions of
\autoref{eq:noslip_bcs} in \autoref{subsec:poiseuille_inout}.

\subsection{Shear flow}\label{subsec:shear}
We use a shear flow boundary condition on all the boundaries so that the
exact solution is $\uu_\mathrm{ref}(\mathbf{x}) = (y,0)$.  With this
boundary condition, the density function $\mathbf{\boldsymbol{\sigma}}$ is only non-zero along
the outer boundary $\Gamma_{0}$.

\subsubsection{$M = 22$ pores}\label{subsec:shear_M22}

\begin{figure}[hbtp!]
  \begin{center}
  (a)~\includegraphics*[height=0.325\linewidth,clip=true]{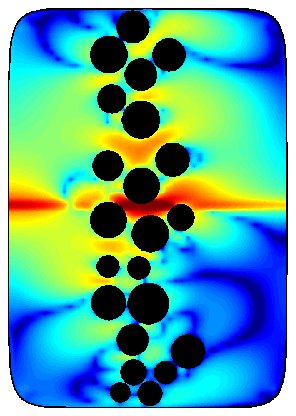}
  \hspace{1em}
  (b)~\includegraphics*[height=0.325\linewidth,clip=true]{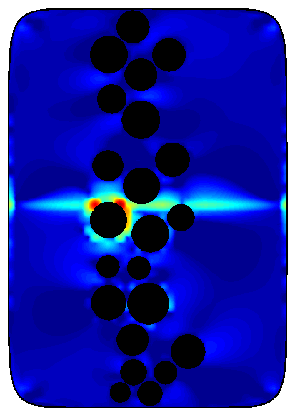}
	  \includegraphics*[height=0.325\linewidth,clip=true]{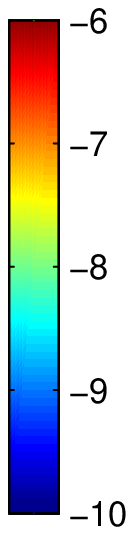}
  \end{center}
  \caption{Logarithm of the relative velocity magnitude error
  $\boldsymbol{\varepsilon}(\xx) = \big| \| \uu(\xx)\| -
  \|\uu_\mathrm{ref}(\xx) \| \big| / \| \uu_\mathrm{ref}(\xx) \|$
  obtained with (a) the block-diagonal preconditioner and (b) the IFMM
  preconditioner, for a shear flow problem with $M = 22$ pores and $N =
  10 \, \mathrm{k}$ unknowns.}
  \label{fig:Stokes_shear_N_9728_errnorm_rel}
\end{figure}

A first small example involves only the first 22 pores (see
\autoref{fig:Stokes_shear_N_9728_errnorm_rel}); the geometry of the
outer boundary $\Gamma_0$ is scaled so that the area of the pore-free
region is the same as the area of the pore-free region in
\autoref{fig:schematic}.  We discretize each pore with $N_{\mathrm{int}}
= 128$ points and the outer boundary with $N_{\mathrm{ext}} = 2048$
points.  The resulting dense system has nearly $10 \, \mathrm{k}$
unknowns.  Three different preconditioning approaches are used in GMRES:
(i)~no preconditioner, (ii)~a block-diagonal preconditioner (with the
diagonal blocks corresponding to the self-interactions of the individual
pores and the outer boundary)\footnote{The diagonal blocks are
assembled, factorized, and inverted exactly.  Alternatively, this could
be accelerated using the FFT since the single-layer potential is nearly
diagonal in Fourier space (see, for example,~\cite{hsi-kop-wen1980},
Theorem 4.1).}, and (iii)~the IFMM.  In the IFMM, ten Chebyshev nodes
are employed in each direction to construct the initial low-rank
approximations ($n = 10$), along with a relative accuracy $\varepsilon =
10^{-7}$ for the low-rank (re)compressions. The preconditioners in cases
(ii) and (iii) are applied from the left as $\mathbf{P}^{-1} \mathbf{A}
\ssigma = \mathbf{P}^{-1} \ff$ where $\mathbf{P}^{-1}$ represents the
preconditioner.

\begin{figure}[hbtp!]
  \begin{center}
  (a)\hspace{-1em}\includegraphics*[width=0.475\linewidth,clip=true]{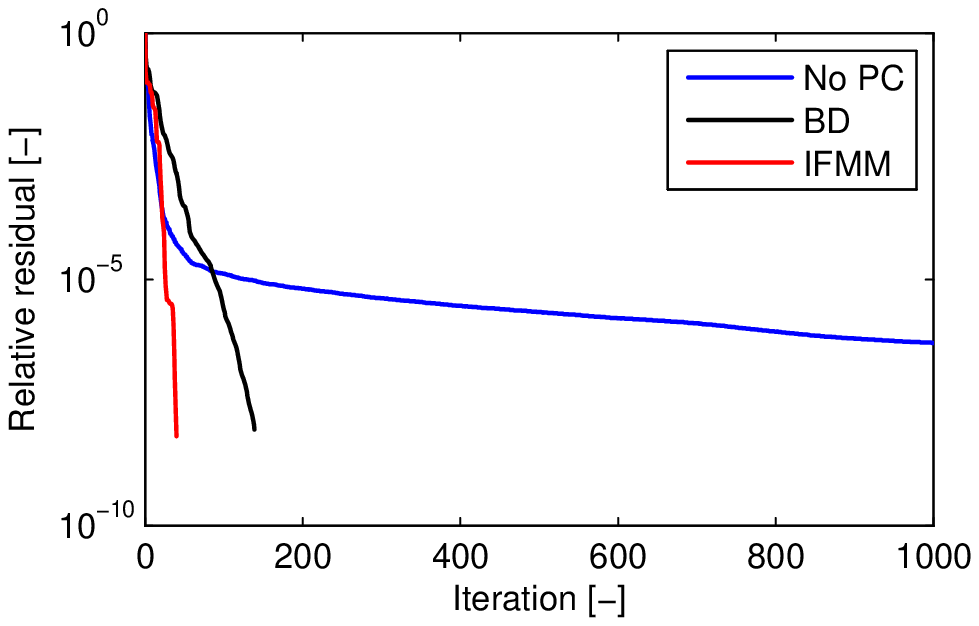}
  (b)\hspace{-1em}\includegraphics*[width=0.475\linewidth,clip=true]{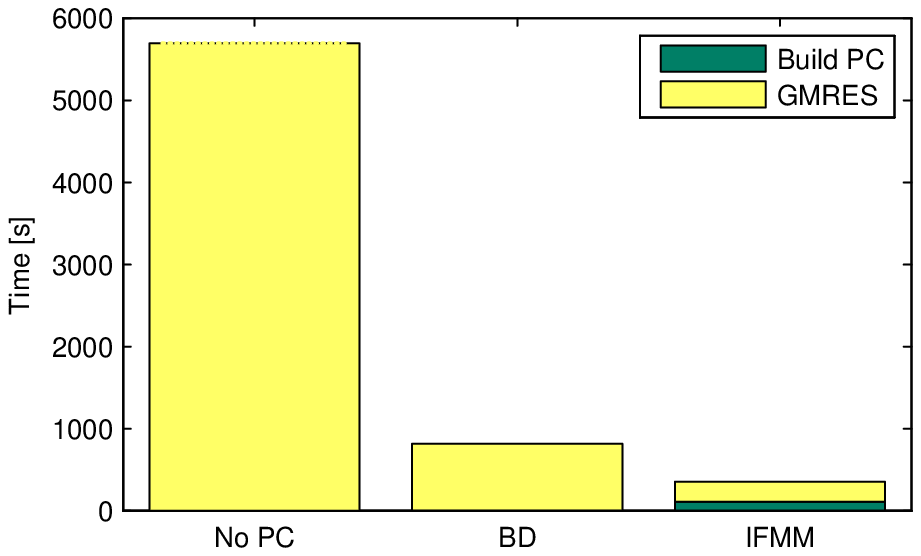}
  \end{center}
  \caption{(a)~Relative residual as a function of the iteration step and
  (b)~total CPU time obtained without a preconditioner (\texttt{No PC}),
  a block-diagonal preconditioner (\texttt{BD}), and the IFMM as the
  preconditioner (\texttt{IFMM}) for a shear flow problem with $M = 22$
  pores and $N = 10 \, \mathrm{k}$ unknowns.  The number of GMRES
  iterations is restricted to 1000.}
  \label{fig:GMRES_Stokes_shear_N_9728}
\end{figure}

\autoref{fig:GMRES_Stokes_shear_N_9728}(a) shows the relative residual
as a function of the number of GMRES iterations and
\autoref{fig:GMRES_Stokes_shear_N_9728}(b) shows the CPU time required
by the preconditioner and GMRES for all three preconditioning
strategies.  Convergence is very slow without a preconditioner and the
relative residual does not drop below the desired tolerance within the
prescribed maximum number of 1000 iterations.  Using the block-diagonal
preconditioner, on the other hand, leads to convergence within 139
iterations. This is a very cheap preconditioning strategy, as the time
needed to construct the preconditioner is negligible. The IFMM
preconditioner is more expensive to compute, but the overall wall time
is significantly reduced because the total number of iterations is
reduced to merely 40.  Note that the relative residual shown in
\autoref{fig:GMRES_Stokes_shear_N_9728}(a) is the left preconditioned
residual $\| \mathbf{P}^{-1}\ff - \mathbf{P}^{-1}\mathbf{A}
\widehat{\ssigma} \|_{2}/ \| \mathbf{P}^{-1} \ff \|_{2}$, where
$\widehat{\ssigma}$ is the estimated solution at each GMRES iteration.
The actual residual $\| \ff - \mathbf{A} \widehat{\ssigma} \|_{2}/ \|
\ff \|_{2}$ is computed as well once the latter drops below the
prescribed tolerance.  \autoref{tbl:GMRES_Stokes_shear_N_9728_tbl}
indicates that the discrepancy between the preconditioned and the actual
relative residual is very small in all cases.

The overall computation times for each of the preconditioning approaches
are presented in \autoref{fig:GMRES_Stokes_shear_N_9728}(b) and
\autoref{tbl:GMRES_Stokes_shear_N_9728_tbl}, where a decomposition is
made into the time required to construct the preconditioner and the
actual iteration time. It is clear that the additional cost for
constructing the IFMM pays off since the total time is more than halved
with respect to the block-diagonal preconditioner.

\begin{table}[htbp]
\footnotesize
\caption{Overview of the results for various preconditioning strategies
for a shear flow problem with $M = 22$ pores and $N = 10 \, \mathrm{k}$
unknowns. The number of GMRES iterations is restricted to 1000.}
\label{tbl:GMRES_Stokes_shear_N_9728_tbl}
\centering
\begin{tabular}{c| c c c c}
	\toprule
    	&	\# iterations	&	Total CPU time [s] & 	$\| \mathbf{P}^{-1}\ff -  \mathbf{P}^{-1}\mathbf{A} \widehat{\ssigma} \|_{2}/ \| \mathbf{P}^{-1} \ff \|_{2}$	  &	$\| \ff -  \mathbf{A} \widehat{\ssigma} \|_{2}/ \| \ff \|_{2}$\\
    	&	&	(Build PC + GMRES)& 	&	\\
	\midrule
	No PC 	& 1000	& $5.69 \times 10^3$ (0 + $5.69 \times 10^3$) 	& $5.16 \times 10^{-7}$ 	& $5.16 \times 10^{-7}$ \\
	BD		& 139	& $8.16 \times 10^2$ ($2.51 \times 10^0$ + $8.13 \times 10^2$)  	& $8.87 \times 10^{-9}$ 	& $8.89 \times 10^{-9}$ \\
	IFMM 	& 40	& $3.52 \times 10^2$ ($1.11 \times 10^2$ + $2.41 \times 10^2$) 	& $6.54 \times 10^{-9}$ 	& $6.54 \times 10^{-9}$ \\
	\bottomrule
\end{tabular}
\end{table}

\begin{figure}[hbtp!]
  \begin{center}
  (a)\hspace{-1em}\includegraphics*[width=0.475\linewidth,clip=true]{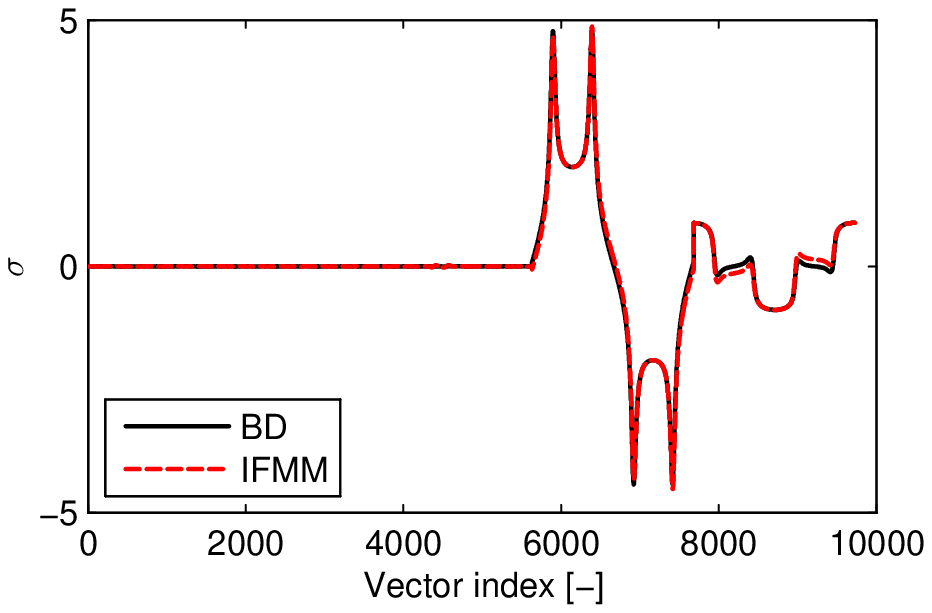}
  (b)\hspace{-1em}\includegraphics*[width=0.475\linewidth,clip=true]{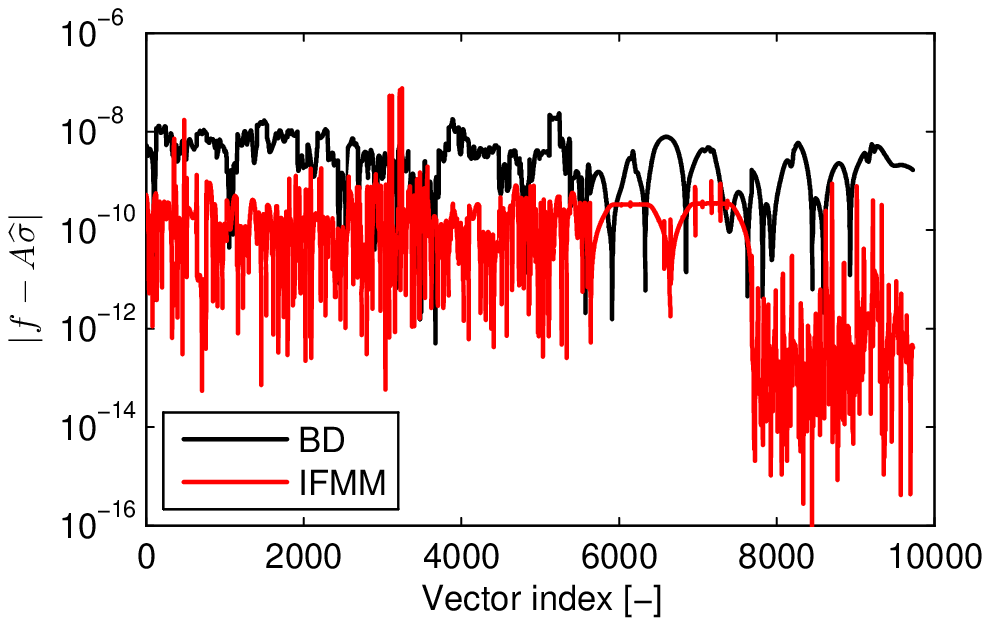}
  \end{center}
  \caption{(a) Density vector $\ssigma$ and (b) absolute value of the
  residual vector $\mathbf{r} = \ff -  \mathbf{A} \ssigma$, obtained
  with a block-diagonal preconditioner (\texttt{BD}), and the IFMM as
  preconditioner (\texttt{IFMM}), for a shear flow problem with $M = 22$
  pores and $N = 10 \, \mathrm{k}$ unknowns.}
  \label{fig:Stokes_shear_N_9728_density_residual}
\end{figure}

The density vector $\ssigma$ and the absolute value of the residual
vector $\mathbf{r} = \ff -  \mathbf{A}\ssigma$ are shown in
\autoref{fig:Stokes_shear_N_9728_density_residual}(a) and
\autoref{fig:Stokes_shear_N_9728_density_residual}(b) for the
block-diagonal and IFMM preconditioner, respectively.  Because of the
boundary condition chosen, the density function is exactly zero on the
boundaries of the pores (the first $2 \times 128 \times 22 = 5632$
degrees of freedom). Generally speaking, the residual is somewhat
smaller for the IFMM than for the block-diagonal preconditioner.

The accuracy of the solutions is further assessed with the logarithm of
the relative error in the magnitude of the velocity field
\begin{align*}
  \boldsymbol{\varepsilon}(\xx) = \frac{\big| \|\uu(\xx)\| -
  \|\uu_\mathrm{ref}(\xx)\| \big|}{\|\uu_\mathrm{ref}(\xx)\|}.
\end{align*}
We plot the error with both the block-diagonal and the IFMM
preconditioner in \autoref{fig:Stokes_shear_N_9728_errnorm_rel}.  We see
that the errors are small in both cases, although the IFMM leads to
errors that are almost two full orders of magnitude smaller.  The
largest errors are observed in both cases along the line $y = 0$ where
the magnitude of the reference velocity field is exactly zero.

%

\subsubsection{$M = 226$ pores}

\begin{figure}[hbtp!]
  \begin{center}
  (a)~\includegraphics*[height=0.275\linewidth,clip=true]{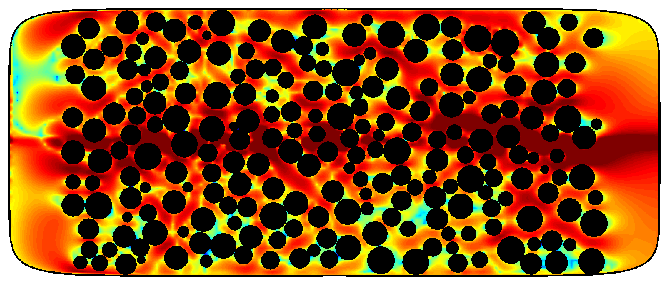}
	  \phantom{\includegraphics*[height=0.275\linewidth,clip=true]{fcolorbar_-10_-6.eps}}\\
  (b)~\includegraphics*[height=0.275\linewidth,clip=true]{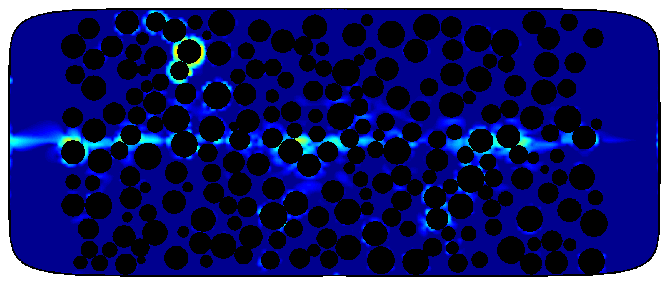}
    \includegraphics*[height=0.275\linewidth,clip=true]{fcolorbar_-10_-6.eps}
  \end{center}
  \caption{Logarithm of the relative velocity magnitude error
  $\boldsymbol{\varepsilon}(\xx) = \big| \|\uu(\xx) \| -
  \|\uu_\mathrm{ref}(\xx) \| \big| / \| \uu_\mathrm{ref}(\xx) \| $
  obtained using (a) the block-diagonal preconditioner and (b) the IFMM
  preconditioner, for a shear flow problem with $M = 226$ pores and $N =
  66 \, \mathrm{k}$ unknowns.}
  \label{fig:Stokes_shear_N_66048_errnorm_rel}
\end{figure}

We now consider the geometry in
\autoref{fig:Stokes_shear_N_66048_errnorm_rel} which contains 226 pores.
We still discretize each pore with $N_{\mathrm{int}}=128$ points, but
the longer outer wall is discretized with $N_{\mathrm{ext}} = 4096$
points, leading to a problem with $66 \, \mathrm{k}$ unknowns. The
block-diagonal preconditioner and the IFMM are used as preconditioners
in GMRES. The number of Chebyshev nodes in the IFMM is first chosen as
$n=10$, but $n=15$ is also investigated since the matrix is expected to
be more ill-conditioned than in the previous example.  Results for
unpreconditioned GMRES are not reported since it requires over 1000
iterations.

\begin{figure}[hbtp!]
  \begin{center}
  (a)\hspace{-1em}\includegraphics*[width=0.45\linewidth,clip=true]{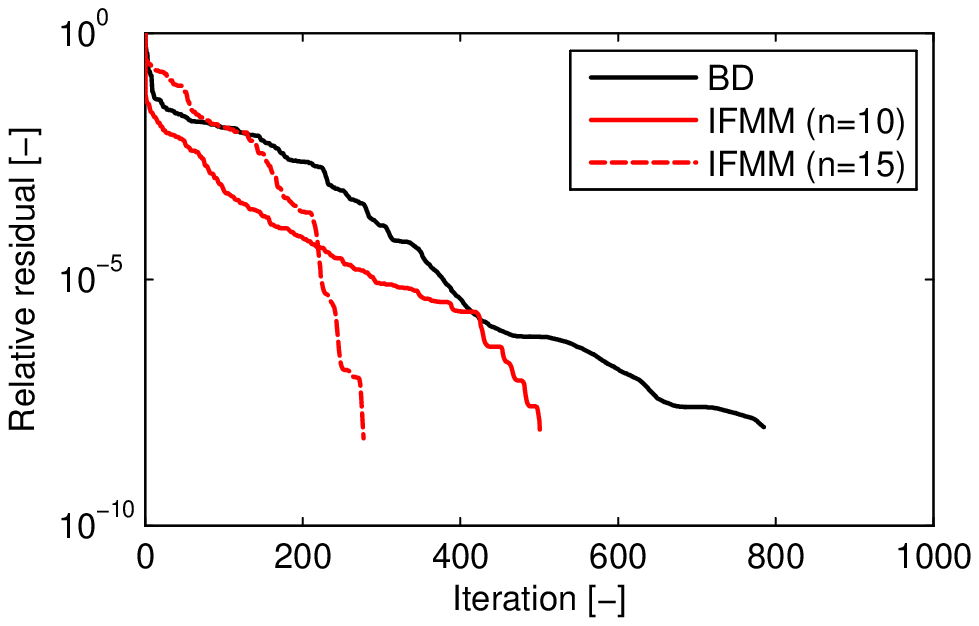}
  (b)\hspace{-1em}\includegraphics*[width=0.45\linewidth,clip=true]{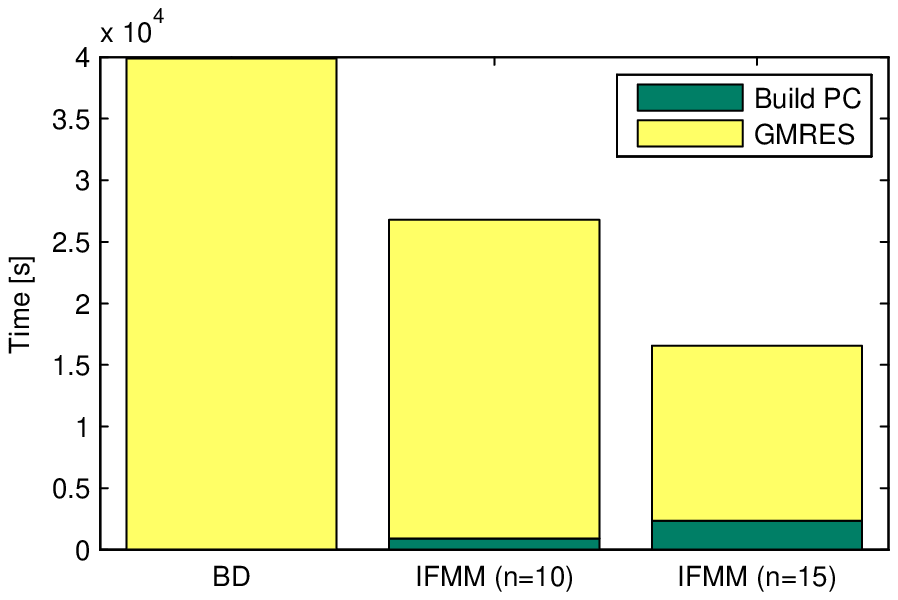}
  \end{center}
  \caption{(a)~Relative residual as a function of the iteration step and
  (b)~total CPU time obtained with a block-diagonal preconditioner
  (\texttt{BD}), and the IFMM as preconditioner (\texttt{IFMM}), for a
  shear flow problem with $M = 226$ pores and $N = 66 \, \mathrm{k}$
  unknowns.}
  \label{fig:GMRES_Stokes_shear_N_66048}
\end{figure}

\autoref{fig:GMRES_Stokes_shear_N_66048}(a) shows the relative residual
as a function of the number of GMRES iterations for all the
preconditioning strategies. GMRES with a block-diagonal preconditioner
requires 785 iterations to converge to the specified residual, which is
significantly larger than the number required in the previous example
(see \autoref{fig:GMRES_Stokes_shear_N_9728}(a)). This demonstrates that
a preconditioner that does not include inter-dependence between pores is
too simplistic to be effective for this complex geometry. Application of
the IFMM, on the other hand, results in a faster convergence of the
residual.  The IFMM preconditioner with $n=10$ requires 501 iterations,
and with $n=15$ requires 277 iterations.

The computation time of each preconditioning approach is depicted in
\autoref{fig:GMRES_Stokes_shear_N_66048}(b).  The most expensive
preconditioner to construct is the IFMM with $n=15$.  However, because
of the reduction in the number of GMRES iterations, this preconditioner
leads to the smallest overall computation time and is more than twice as
fast as the block-diagonal preconditioner.


%

Apart from achieving a speedup, application of the IFMM also leads to a
more accurate solution than the block-diagonal preconditioner. This is
highlighted in the two last columns of
\autoref{tbl:GMRES_Stokes_shear_N_66048_tbl}: although the
preconditioned residual $\| \mathbf{P}^{-1}\ff -
\mathbf{P}^{-1}\mathbf{A} \widehat{\ssigma} \|_{2}/ \| \mathbf{P}^{-1}
\ff \|_{2}$ has reached the desired tolerance for the block-diagonal
preconditioner, the actual residual $\| \ff -  \mathbf{A}
\widehat{\ssigma} \|_{2}/ \| \ff \|_{2}$ is still several orders of
magnitude larger.

\begin{table}[htbp!]
\footnotesize
\caption{Overview of the results for various preconditioning strategies
for a shear flow problem with $M = 226$ pores and $N = 66 \, \mathrm{k}$
unknowns.}
\label{tbl:GMRES_Stokes_shear_N_66048_tbl}
\centering
\begin{tabular}{c| c c c c}
	\toprule
    	&	\# iterations	&	Total CPU time [s] & 	$\| \mathbf{P}^{-1}\ff -  \mathbf{P}^{-1}\mathbf{A} \widehat{\ssigma} \|_{2}/ \| \mathbf{P}^{-1} \ff \|_{2}$	  &	$\| \ff -  \mathbf{A} \widehat{\ssigma} \|_{2}/ \| \ff \|_{2}$\\
    	&	&	(Build PC + GMRES)& 	&	\\
	\midrule
	BD				& 785	& $3.99 \times 10^4$ ($1.78 \times 10^1$ + $3.99 \times 10^4$)  		& $1.00 \times 10^{-8}$ 	& $3.67 \times 10^{-2}$ \\
	IFMM ($n=10$) 	& 501	& $2.68 \times 10^4$ ($8.93 \times 10^2$ + $2.59 \times 10^4$) 		& $8.94 \times 10^{-9}$ 	& $8.94 \times 10^{-9}$ \\
	IFMM ($n=15$)	& 277	& $1.65 \times 10^4$ ($2.34 \times 10^3$ + $1.42 \times 10^4$) 		& $5.94 \times 10^{-9}$ 	& $5.94 \times 10^{-9}$ \\
	\bottomrule
\end{tabular}
\end{table}

The accuracy of the solution is further  assessed in
\autoref{fig:Stokes_shear_N_66048_errnorm_rel}, showing the logarithm
of the relative velocity magnitude error
$\boldsymbol{\varepsilon}(\xx)$.  The largest errors appear once again
along the line $y=0$. Furthermore, the errors obtained with the IFMM
are almost four orders of magnitude smaller than with the
block-diagonal preconditioner. The results presented in
\autoref{fig:GMRES_Stokes_shear_N_66048} and
\ref{fig:Stokes_shear_N_66048_errnorm_rel} clearly demonstrate the
efficiency and accuracy of the IFMM for the challenging problem under
concern.

%

\subsubsection{$M = 826$ pores}
Finally, the full geometry with 826 pores depicted in
\autoref{fig:schematic} is considered.  We discretize each pore with
$N_{\mathrm{int}} = 128$ points and the outer boundary with
$N_{\mathrm{ext}} = 16384$ points, leading to a problem with $244 \,
\mathrm{k}$ unknowns. The number of Chebyshev nodes in the IFMM is chosen as $n=15$, while the maximum number of GMRES iterations is increased to 2000.

\begin{figure}[hbtp!]
  \begin{center}
  \includegraphics*[width=0.95\linewidth,clip=true]{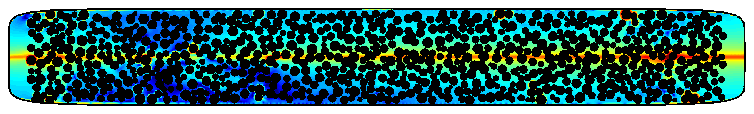}
  \includegraphics*[height=0.125\linewidth,clip=true]{fcolorbar_-10_-6.eps}
  \end{center}
  \caption{Logarithm of the relative velocity magnitude error
  $\boldsymbol{\varepsilon}(\xx) = \big| \|\uu(\xx) \| -
  \|\uu_\mathrm{ref}(\xx) \| \big| / \| \uu_\mathrm{ref}(\xx) \| $
  obtained using the IFMM
  preconditioner, for a shear flow problem with $M = 826$ pores and $N =
  244 \, \mathrm{k}$ unknowns.}
  \label{fig:Stokes_shear_N_244224_errnorm_rel}
\end{figure}

\autoref{tbl:GMRES_Stokes_shear_N_244224_tbl} indicates that the block-diagonal preconditioner does not reach convergence within 2000 iterations, and that the actual residual remains very large. On the other hand, convergence is obtained after 1489 iterations with the IFMM (with a small actual residual), confirming its effectiveness even for large problems with a very complicated geometry. \autoref{fig:Stokes_shear_N_244224_errnorm_rel} depicts the logarithm of the relative velocity magnitude error $\boldsymbol{\varepsilon}(\xx)$, as obtained with the IFMM preconditioner. The errors are larger than in the previous test cases (compare with \autoref{fig:Stokes_shear_N_9728_errnorm_rel}(b) and \autoref{fig:Stokes_shear_N_66048_errnorm_rel}(b)), but are still reasonably small.

\begin{table}[htbp!]
\footnotesize
\caption{Overview of the results for various preconditioning strategies
for a shear flow problem with $M = 826$ pores and $N = 244 \, \mathrm{k}$
unknowns. The number of GMRES iterations is restricted to 2000.}
\label{tbl:GMRES_Stokes_shear_N_244224_tbl}
\centering
\begin{tabular}{c| c c c c}
	\toprule
    	&	\# iterations	&	Total CPU time [s] & 	$\| \mathbf{P}^{-1}\ff -  \mathbf{P}^{-1}\mathbf{A} \widehat{\ssigma} \|_{2}/ \| \mathbf{P}^{-1} \ff \|_{2}$	  &	$\| \ff -  \mathbf{A} \widehat{\ssigma} \|_{2}/ \| \ff \|_{2}$\\
    	&	&	(Build PC + GMRES)& 	&	\\
	\midrule
	BD				& 	2000 &  $3.04 \times 10^5$ ($3.57 \times 10^2$ + $3.04 \times 10^5$)  		& $7.44 \times 10^{-6}$ 	& $1.82 \times 10^{-1}$ \\
	IFMM ($n=15$)	& 	1489 &  $2.70 \times 10^5$ ($8.83 \times 10^3$ + $2.61 \times 10^5$) 		& $9.98 \times 10^{-9}$ 	& $1.15 \times 10^{-6}$ \\
	\bottomrule
\end{tabular}
\end{table}
%

\subsection{A linear combination of Stokeslets and rotlets}
\label{subsec:analytical}

The examples considered in \autoref{subsec:shear} are repeated with the
same number of discretization points ($N_{\mathrm{int}} = 128$ points
per pore and $N_{\mathrm{ext}} = 2048, 4096, 16384$ for the outer
boundary), but with a boundary condition that results in a non-zero
density vector $\ssigma$ on the boundaries of the inner pores.  The
boundary condition is given by a linear combination of Stokeslets and
rotlets
\begin{align}
  \ff(\xx) =\sum_{k=0}^{M}\left( -\log \| \xx - \mathbf{c}_k \| \mathbf{I} +
  \frac{\left(\xx - \mathbf{c}_k \right) \otimes \left(\xx -
  \mathbf{c}_k \right)}{\| \xx - \mathbf{c}_k \|^{2}} \right)
  \boldsymbol{\lambda}_k + \sum_{k=0}^{M}\frac{\left(\xx - \mathbf{c}_k
  \right)^{\perp}}{\| \xx - \mathbf{c}_k \|^{2}} \mu_k, \quad \xx \in
  \Gamma
  \label{eq:analyticalBC}
\end{align}
where $\mathbf{c}_0$ is an arbitrary point located outside the domain
$\Omega_0$ and $\mathbf{c}_k$ ($k=1 \ldots M$) are arbitrary points
inside each of the pores $\Omega_{k}$. The vectors
$\boldsymbol{\lambda}_k$ and scalars $\mu_k$ are randomly chosen, and
$(x,y)^{\perp} = (y,-x)$.  The exact velocity field
$\uu_\mathrm{ref}(\xx)$ at any point $\xx \in \Omega$ is given by
\autoref{eq:analyticalBC} as well.


\subsubsection{$M = 22$ pores}
We again consider the first 22 pores.  The GMRES convergence results
are very similar to the shear flow example: no convergence is achieved
in 1000 iterations without preconditioning, while 143 and only 54
iterations are needed with the block-diagonal and the IFMM
preconditioners, respectively. The results are summarized in
\autoref{tbl:GMRES_Stokes_analytical_N_9728_tbl}.  We see that the IFMM
preconditioner is nearly twice as fast as the block-diagonal
preconditioner.

\begin{table}[htbp!]
\footnotesize
\caption{Overview of the results for various preconditioning strategies
for a Stokeslets/rotlets flow problem with $M = 22$ pores and $N = 10 \,
\mathrm{k}$ unknowns. The number of GMRES iterations is restricted to
1000.}
\label{tbl:GMRES_Stokes_analytical_N_9728_tbl}
\centering
\begin{tabular}{c| c c c c}
	\toprule
    	&	\# iterations	&	Total CPU time [s] & 	$\| \mathbf{P}^{-1}\ff -  \mathbf{P}^{-1}\mathbf{A} \widehat{\ssigma} \|_{2}/ \| \mathbf{P}^{-1} \ff \|_{2}$	  &	$\| \ff -  \mathbf{A} \widehat{\ssigma} \|_{2}/ \| \ff \|_{2}$\\
    	&	&	(Build PC + GMRES)& 	&	\\
	\midrule
	No PC 	& 1000	& $5.65 \times 10^3$ (0 + $5.65 \times 10^3$) 	& $2.39 \times 10^{-5}$ 	& $2.39 \times 10^{-5}$ \\
	BD		& 143	& $8.39 \times 10^2$ ($2.44 \times 10^0$ + $8.36 \times 10^2$)  	& $8.61 \times 10^{-9}$ 	& $8.65 \times 10^{-9}$ \\
	IFMM 	& 55	& $4.37 \times 10^2$ ($1.20 \times 10^2$ + $3.18 \times 10^2$) 	& $3.79 \times 10^{-9}$ 	& $3.79 \times 10^{-9}$ \\
	\bottomrule
\end{tabular}
\end{table}

\begin{figure}[hbtp!]
  \begin{center}
  (a)\hspace{-0em}\includegraphics*[height=0.325\linewidth,clip=true]{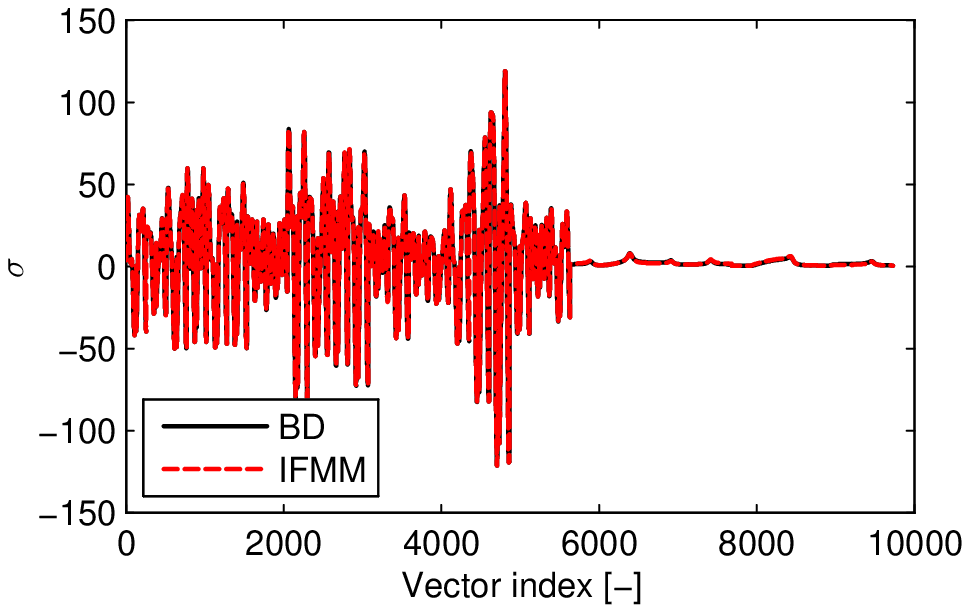}\hspace{1em}
  (b)\hspace{-0em}\includegraphics*[height=0.325\linewidth,clip=true]{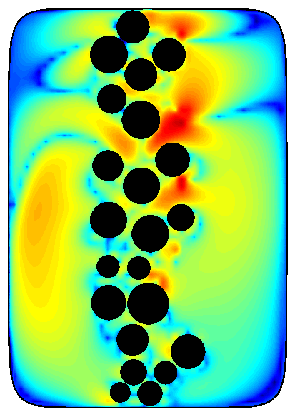}
  	  \includegraphics*[height=0.325\linewidth,clip=true]{fcolorbar_-10_-6.eps}
  \end{center}
  \caption{(a) Density vector $\ssigma$ obtained with a block-diagonal
  preconditioner (\texttt{BD}) and the IFMM as preconditioner
  (\texttt{IFMM}), and (b) logarithm of the relative velocity magnitude
  error $\boldsymbol{\varepsilon}(\xx) = \big| \|\uu(\xx)\| - \|
  \uu_\mathrm{ref}(\xx) \| \big| / \| \uu_\mathrm{ref}(\xx) \|$ obtained
  with the IFMM preconditioner, for a linear combination of Stokeslets
  and rotlets with $M = 22$ pores and $N = 10 \, \mathrm{k}$ unknowns.}
  \label{fig:Stokes_analytical_N_9728_density_error}
\end{figure}

With this boundary condition, the density vector $\ssigma$ is non-zero
on the boundary of each inner pore
(\autoref{fig:Stokes_analytical_N_9728_density_error}(a)). A good
correspondence between the results of both preconditioners is found.
\autoref{fig:Stokes_analytical_N_9728_density_error}(b) shows the
logarithm of the relative velocity magnitude error
$\boldsymbol{\varepsilon}(\xx)$ obtained with the IFMM.  As in the shear flow
example, about six digits of accuracy is achieved.

%

\subsubsection{$M = 226$ pores}
We next increase the problem size to the first 226 pores and summarize
the results in \autoref{tbl:GMRES_Stokes_analytical_N_66048_tbl}. GMRES
with the block-diagonal preconditioner needs over 700 iterations to
converge, while application of the IFMM reduces this number to 506 for
$n=10$ and 258 for $n=15$. As we saw for the shear boundary condition,
using $n=15$ Chebyshev nodes results in the smallest overall computation
time, and this is almost 2.5 times faster than the block-diagonal
preconditioner.

\begin{table}[htbp!]
\footnotesize
\caption{Overview of the results for various preconditioning strategies
for a Stokeslets/rotlets flow problem with $M = 226$ pores and $N = 66
\, \mathrm{k}$ unknowns.}
\label{tbl:GMRES_Stokes_analytical_N_66048_tbl}
\centering
\begin{tabular}{c| c c c c}
	\toprule
    	&	\# iterations	&	Total CPU time [s] & 	$\| \mathbf{P}^{-1}\ff -  \mathbf{P}^{-1}\mathbf{A} \widehat{\ssigma} \|_{2}/ \| \mathbf{P}^{-1} \ff \|_{2}$	  &	$\| \ff -  \mathbf{A} \widehat{\ssigma} \|_{2}/ \| \ff \|_{2}$\\
    	&	&	(Build PC + GMRES)& 	&	\\
	\midrule
	BD				& 715	& $3.63 \times 10^4$ ($1.81 \times 10^1$ + $3.63 \times 10^4$)  		& $9.94 \times 10^{-9}$ 	& $6.07 \times 10^{-3}$ \\
	IFMM ($n=10$)	& 506	& $2.68 \times 10^4$ ($8.94 \times 10^2$ + $2.59 \times 10^4$) 		& $8.33 \times 10^{-9}$ 	& $8.33 \times 10^{-9}$ \\
	IFMM ($n=15$)	& 258	& $1.57 \times 10^4$ ($2.35 \times 10^3$ + $1.33 \times 10^4$) 		& $7.33 \times 10^{-9}$ 	& $7.33 \times 10^{-9}$ \\
	\bottomrule
\end{tabular}
\end{table}

\begin{figure}[hbtp!]
  \begin{center}
  (a)~\includegraphics*[height=0.275\linewidth,clip=true]{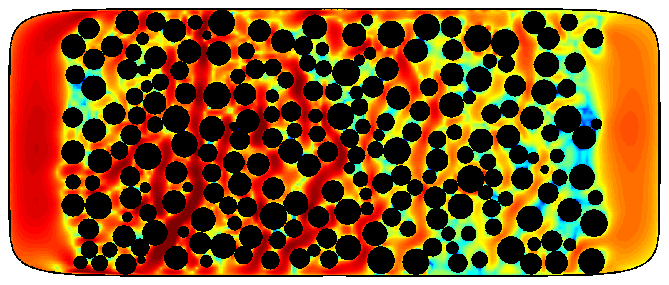}
  	  \phantom{\includegraphics*[height=0.275\linewidth,clip=true]{fcolorbar_-10_-6.eps}}\\
  (b)~\includegraphics*[height=0.275\linewidth,clip=true]{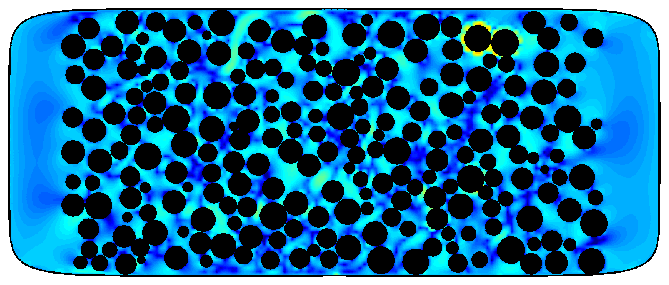}
	  \includegraphics*[height=0.275\linewidth,clip=true]{fcolorbar_-10_-6.eps}
  \end{center}
  \caption{Logarithm of the relative velocity magnitude error
  $\boldsymbol{\varepsilon}(\xx) = \big| \|\uu(\xx)\| - \|
  \uu_\mathrm{ref}(\xx) \| \big| / \| \uu_\mathrm{ref}(\xx) \| $
  obtained using (a) the block-diagonal preconditioner and (b) the IFMM
  preconditioner, for a Stokeslets/rotlets flow problem with $M = 226$
  pores and $N = 66 \, \mathrm{k}$ unknowns.}
  \label{fig:Stokes_analytical_N_66048_errnorm_rel}
\end{figure}

Similar to the shear flow example, there is a large discrepancy between
the preconditioned and the actual residual for the block-diagonal
preconditioner, while this is not the case for the IFMM (see
\autoref{tbl:GMRES_Stokes_analytical_N_66048_tbl}). The errors inside
the domain are consequently several orders of magnitude smaller when
using the IFMM preconditioner as is illustrated in
\autoref{fig:Stokes_analytical_N_66048_errnorm_rel}.
%

\subsubsection{$M = 826$ pores}

\begin{figure}[hbtp!]
  \begin{center}
  \includegraphics*[width=0.95\linewidth,clip=true]{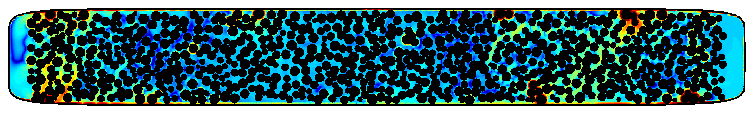}
  \includegraphics*[height=0.125\linewidth,clip=true]{fcolorbar_-10_-6.eps}
  \end{center}
  \caption{Logarithm of the relative velocity magnitude error
  $\boldsymbol{\varepsilon}(\xx) = \big| \|\uu(\xx) \| -
  \|\uu_\mathrm{ref}(\xx) \| \big| / \| \uu_\mathrm{ref}(\xx) \| $
  obtained using the IFMM preconditioner, for a Stokeslets/rotlets flow problem with $M = 826$ pores and $N =
  244 \, \mathrm{k}$ unknowns.}
  \label{fig:Stokes_analytical_N_244224_errnorm_rel}
\end{figure}

Finally, all $826$ pores are considered.
\autoref{tbl:GMRES_Stokes_analytical_N_244224_tbl} indicates that no
convergence is achieved with the block-diagonal preconditioner (even
after 2000 iterations), while 1283 iterations are needed with the IFMM.
\autoref{fig:Stokes_analytical_N_244224_errnorm_rel} confirms the
accuracy of the solution obtained with the IFMM, although some localized
spots with larger errors are observed.

\begin{table}[htbp!]
\footnotesize
\caption{Overview of the results for various preconditioning strategies
for a Stokeslets/rotlets flow problem with $M = 826$ pores and $N = 244
\, \mathrm{k}$ unknowns. The number of GMRES iterations is restricted to 2000.}
\label{tbl:GMRES_Stokes_analytical_N_244224_tbl}
\centering
\begin{tabular}{c| c c c c}
	\toprule
    	&	\# iterations	&	Total CPU time [s] & 	$\| \mathbf{P}^{-1}\ff -  \mathbf{P}^{-1}\mathbf{A} \widehat{\ssigma} \|_{2}/ \| \mathbf{P}^{-1} \ff \|_{2}$	  &	$\| \ff -  \mathbf{A} \widehat{\ssigma} \|_{2}/ \| \ff \|_{2}$\\
    	&	&	(Build PC + GMRES)& 	&	\\
	\midrule
	BD				& 2000 &  $3.04 \times 10^5$ ($3.61 \times 10^2$ + $3.04 \times 10^5$)  	& $1.12 \times 10^{-7}$ 	& $7.23 \times 10^{-2}$ \\
	IFMM ($n=15$)	& 1283 &  $2.24 \times 10^5$ ($8.82 \times 10^3$ + $2.16 \times 10^5$) 		& $1.00 \times 10^{-8}$ 	& $1.34 \times 10^{-7}$ \\
	\bottomrule
\end{tabular}
\end{table}
%

\subsection{Pipe flow}
\label{subsec:poiseuille_inout}

This final section considers the more physical boundary conditions
defined in \autoref{eq:noslip_bcs} and illustrated in
\autoref{fig:schematic}.  The velocity is parabolic at the intake and
outtake of the outer boundary $\Gamma_0$, and zero on each pore.  The
same number of discretization points as in \autoref{subsec:shear} and
\ref{subsec:analytical} is used.

\subsubsection{$M = 22$ pores}
\begin{figure}[hbtp!]
  \begin{center}
  (a)\hspace{-1em}\includegraphics*[width=0.475\linewidth,clip=true]{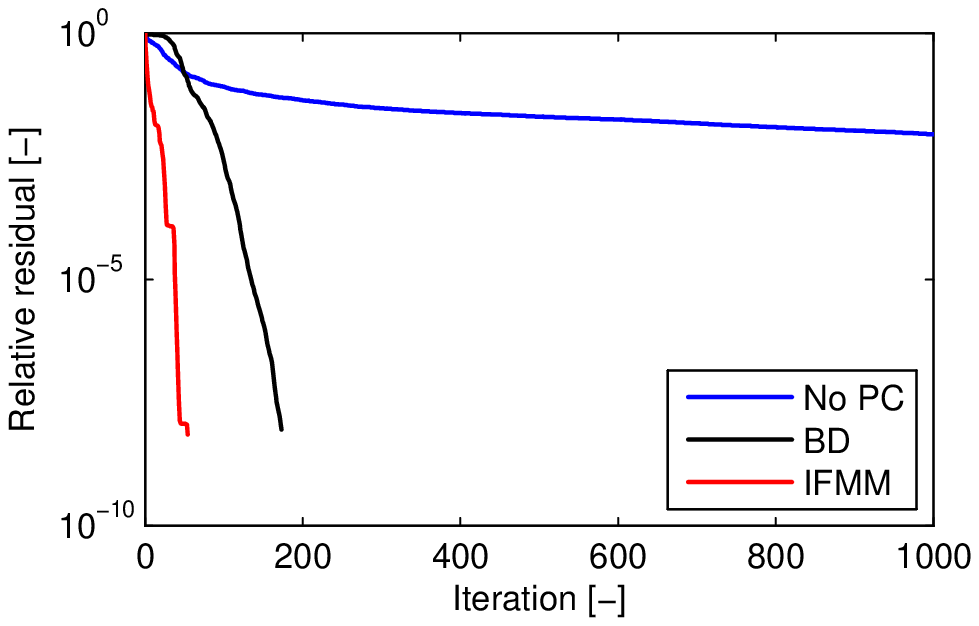}
  (b)\hspace{-1em}\includegraphics*[width=0.475\linewidth,clip=true]{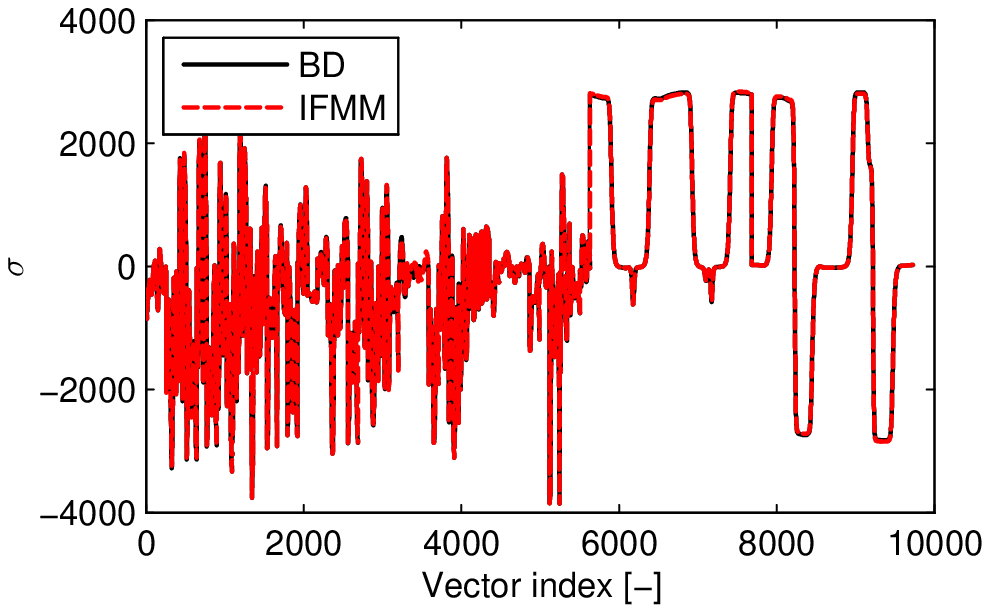}
  \end{center}
  \caption{(a)~Relative residual as a function of the iteration step and
  (b)~density vector $\ssigma$, obtained without preconditioner
  (\texttt{No PC}), a block-diagonal preconditioner (\texttt{BD}), and
  the IFMM as preconditioner (\texttt{IFMM}) for a pipe flow
  problem with $M = 22$ pores and $N = 10 \, \mathrm{k}$ unknowns.}
  \label{fig:GMRES_Stokes_Poiseuille_inout_N_9728}
\end{figure}

\autoref{fig:GMRES_Stokes_Poiseuille_inout_N_9728}(a) shows the relative
residual as a function of the number of GMRES iterations for the case
involving 22 pores. The convergence behavior of the preconditioners is
similar to the previous two examples.  However, without preconditioning,
the convergence is notably slower than before.  The results are
summarized in \autoref{tbl:GMRES_Stokes_Poiseuille_inout_N_9728_tbl}.
The density vector $\ssigma$ obtained with the block-diagonal and IFMM
preconditioners is depicted in
\autoref{fig:GMRES_Stokes_Poiseuille_inout_N_9728}(b), illustrating a
good match between both approaches.

\begin{figure}[hbtp!]
  \begin{center}
   \includegraphics*[scale=0.4,clip=true,bb=400 40 880 530]{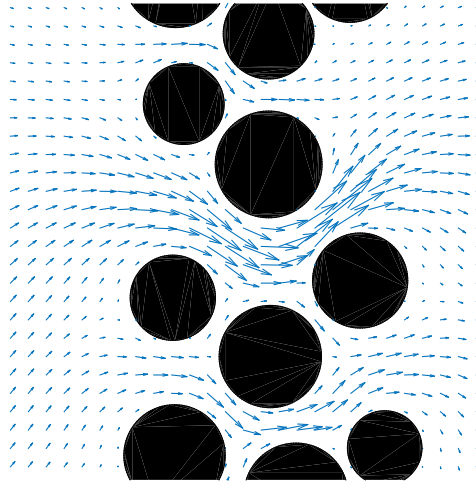}
  \end{center}
  \caption{Quiver plot of the velocity field  $\uu(\xx)$ obtained when
  using the IFMM preconditioner for a pipe flow problem with $M = 22$
  pores and $N = 10 \, \mathrm{k}$ unknowns.  The velocity field is
  plotted in a small window in the middle of the porous region.}
  \label{fig:Stokes_Poiseuille_inout_N_9728_velocityfield}
\end{figure}

\begin{table}[htbp!]
\footnotesize
\caption{Overview of the results for various preconditioning strategies,
for a pipe flow problem with $M = 22$ pores and $N = 10 \, \mathrm{k}$
unknowns. The number of GMRES iterations is restricted to 1000.}
\label{tbl:GMRES_Stokes_Poiseuille_inout_N_9728_tbl}
\centering
\begin{tabular}{c| c c c c}
	\toprule
    	&	\# iterations	&	Total CPU time [s] & 	$\| \mathbf{P}^{-1}\ff -  \mathbf{P}^{-1}\mathbf{A} \widehat{\ssigma} \|_{2}/ \| \mathbf{P}^{-1} \ff \|_{2}$	  &	$\| \ff -  \mathbf{A} \widehat{\ssigma} \|_{2}/ \| \ff \|_{2}$\\
    	&	&	(Build PC + GMRES)& 	&	\\
	\midrule
	No PC 	& 1000	& $5.66 \times 10^3$ (0 + $5.66 \times 10^3$) 	& $8.81 \times 10^{-3}$ 	& $8.81 \times 10^{-3}$ \\
	BD		& 173	& $3.61 \times 10^2$ ($2.44 \times 10^0$ + $3.59 \times 10^2$)  	& $8.95 \times 10^{-9}$ 	& $8.98 \times 10^{-9}$ \\
	IFMM 	& 54	& $9.80 \times 10^1$ ($1.72 \times 10^1$ + $8.08 \times 10^1$) 		& $7.07 \times 10^{-9}$ 	& $7.07 \times 10^{-9}$ \\
	\bottomrule
\end{tabular}
\end{table}
\autoref{fig:Stokes_Poiseuille_inout_N_9728_velocityfield} shows a
quiver plot of the velocity field $\uu(\xx)$ obtained with the
block-diagonal and the IFMM preconditioner.  The parabolic velocity
profile imposed on the left side of the outer boundary is strongly
affected by the presence of the pores, and large velocities are observed
in the areas with the largest gaps between the pores.



%

\subsubsection{$M = 226$ pores}\label{subsubsec:Poisueille_226}

\begin{figure}[hbtp!]
  \begin{center}
  \includegraphics*[width=0.475\linewidth,clip=true]{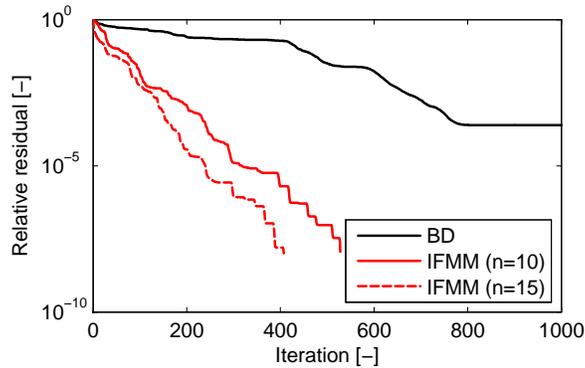}
  \end{center}
  \caption{Relative residual as a function of the iteration step for
  a block-diagonal preconditioner (\texttt{BD}) and the IFMM
  preconditioner (\texttt{IFMM}) for a pipe flow with $M = 226$ pores and $N = 66 \,
  \mathrm{k}$ unknowns. The number of GMRES iterations is restricted to
  1000.}
  \label{fig:GMRES_Stokes_Poiseuille_inout_N_66048}
\end{figure}

For all the examples considered so far, the block-diagonal
preconditioner did not perform that badly; application of the IFMM has
only lead to rather moderate speedups (up to a factor  of 2.5). The
pipe flow problem involving 226 pores turns out to be much harder for
the block-diagonal preconditioner, however, as is indicated in
\autoref{fig:GMRES_Stokes_Poiseuille_inout_N_66048} and
\autoref{tbl:GMRES_Stokes_Poiseuille_inout_N_66048_tbl}. Block-diagonal
preconditioned GMRES is unable to converge within 1000 iterations and
the residual is still more than four orders of magnitude above the
prescribed tolerance of $10^{-8}$. Application of the IFMM, on the other
hand, guarantees converge in 528 iterations for $n=10$ and slightly more
than 400 iterations for $n=15$, hence confirming the effectiveness of
the IFMM.

\begin{table}[htbp!]
\footnotesize
\caption{Overview of the results for various preconditioning strategies
for a pipe flow problem with $M = 226$ pores and $N = 66 \, \mathrm{k}$ unknowns. The
number of GMRES iterations is restricted to 1000.}
\label{tbl:GMRES_Stokes_Poiseuille_inout_N_66048_tbl}
\centering
\begin{tabular}{c| c c c c}
	\toprule
    	&	\# iterations	&	Total CPU time [s] & 	$\| \mathbf{P}^{-1}\ff -  \mathbf{P}^{-1}\mathbf{A} \widehat{\ssigma} \|_{2}/ \| \mathbf{P}^{-1} \ff \|_{2}$	  &	$\| \ff -  \mathbf{A} \widehat{\ssigma} \|_{2}/ \| \ff \|_{2}$\\
    	&	&	(Build PC + GMRES)& 	&	\\
	\midrule
	BD				& 1000	& $4.42 \times 10^4$ ($1.81 \times 10^1$ + $4.42 \times 10^4$)  		& $2.45 \times 10^{-4}$ 	& $9.31 \times 10^{-1}$ \\
	IFMM ($n=10$)	& 528	& $1.30 \times 10^4$ ($5.63 \times 10^2$ + $1.24 \times 10^4$) 		& $9.97 \times 10^{-9}$ 	& $4.64 \times 10^{-6}$ \\
	IFMM ($n=15$)	& 407	& $1.12 \times 10^4$ ($1.81 \times 10^3$ + $9.38 \times 10^3$) 		& $9.94 \times 10^{-9}$ 	& $2.07 \times 10^{-6}$ \\
	\bottomrule
\end{tabular}
\end{table}

\begin{figure}[hbtp!]
  \begin{center}
   \includegraphics*[scale=0.4,clip=true,bb=400 40 880 530]{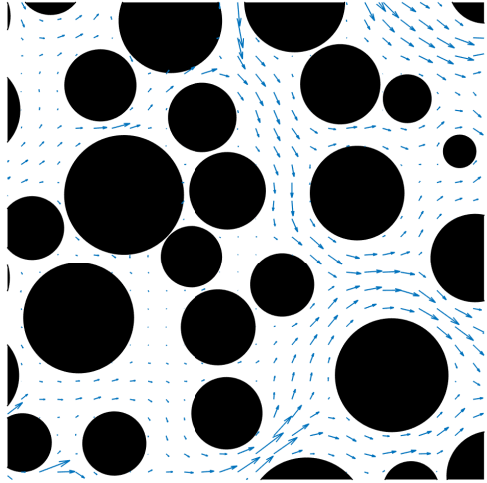}
  \end{center}
  \caption{Quiver plot of the velocity field  $\uu(\xx)$ obtained when
  using the IFMM preconditioner for a pipe flow problem with $M = 226$ pores and $N =
  66 \, \mathrm{k}$ unknowns.  The velocity field is plotted in a small
  window in the middle of the porous region.}
  \label{fig:Stokes_Poiseuille_inout_N_66048_velocityfield}
\end{figure}

\autoref{fig:Stokes_Poiseuille_inout_N_66048_velocityfield} shows a
quiver plot of the velocity field $\uu(\xx)$ obtained with the IFMM
preconditioner. Even though the governing equations are linear, the
flow is very complex because of the geometry.  With this velocity
field, particle tracking can be performed, and this can be used to find
structure within the flow~\cite{hal2001}.  However, this is only
practical if the linear system \autoref{eq:dense_matrix} can be solved
with a specified tolerance with a reasonable amount of computation
time.  With our proposed IFMM preconditioner, this is exactly the case.

\subsubsection{$M = 826$ pores}

Finally, the problem depicted in \autoref{fig:schematic} is solved. The
block diagonal preconditioner is once more unable to reach convergence
and the residual remains large after 2000 iterations, as is summarized
in \autoref{tbl:GMRES_Stokes_Poiseuille_inout_N_244224_tbl}. GMRES would
probably need a few more thousand iterations to converge with this
preconditioner. The IFMM preconditioner needs 1645 iterations to
converge, leading to a reasonably small actual residual.

\begin{table}[htbp!]
\footnotesize
\caption{Overview of the results for various preconditioning strategies
for a pipe flow problem with $M = 826$ pores and $N = 244 \, \mathrm{k}$ unknowns. The
number of GMRES iterations is restricted to 2000.}
\label{tbl:GMRES_Stokes_Poiseuille_inout_N_244224_tbl}
\centering
\begin{tabular}{c| c c c c}
	\toprule
    	&	\# iterations	&	Total CPU time [s] & 	$\| \mathbf{P}^{-1}\ff -  \mathbf{P}^{-1}\mathbf{A} \widehat{\ssigma} \|_{2}/ \| \mathbf{P}^{-1} \ff \|_{2}$	  &	$\| \ff -  \mathbf{A} \widehat{\ssigma} \|_{2}/ \| \ff \|_{2}$\\
    	&	&	(Build PC + GMRES)& 	&	\\
	\midrule
	BD				& 2000	&  $3.04 \times 10^5$ ($3.61 \times 10^2$ + $3.04 \times 10^5$)	& $9.78 \times 10^{-3}$ 	& $9.48 \times 10^{-1}$ \\
	IFMM ($n=15$)	& 1645	&  $2.96 \times 10^5$ ($8.84 \times 10^3$ + $2.87 \times 10^5$) 		& $1.00 \times 10^{-8}$ 	& $2.67 \times 10^{-5}$ \\
	\bottomrule
\end{tabular}
\end{table}

\begin{figure}[hbtp!]
  \begin{center}
   \includegraphics*[scale=0.4,clip=true,bb=400 40 880 530]{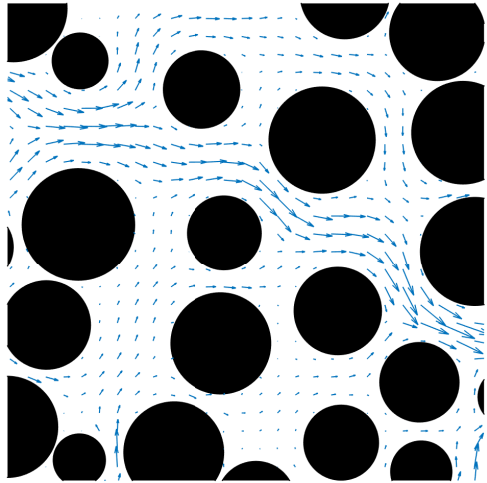}
  \end{center}
  \caption{Quiver plot of the velocity field  $\uu(\xx)$ obtained when
  using the IFMM preconditioner for a pipe flow problem with $M = 826$
  pores and $N = 244 \, \mathrm{k}$ unknowns.  The velocity field is
  plotted in a small window in the middle of the porous region.}
  \label{fig:Stokes_Poiseuille_inout_N_244224_velocityfield}
\end{figure}

\autoref{fig:Stokes_Poiseuille_inout_N_244224_velocityfield} shows a
quiver plot of the velocity field $\uu(\xx)$ obtained with the IFMM
preconditioner. A very complex flow is resolved in this case as well.

%

\section{Conclusions}\label{sec:Conclusions}

Discretization of a first-kind boundary integral equation, representing
a 2D incompressible Stokes flow in a porous medium, leads to a dense
linear system of equations that is computationally expensive to solve.
In particular, a large number of iterations is often required if an
iterative Krylov methods such as GMRES is employed. In this paper, the
IFMM has been presented as an efficient preconditioner that is capable
of significantly reducing the number of iterations and the overall
computation cost. 

The IFMM is in essence an inexact fast direct solver for dense
$\mH^{2}$-matrices with a linear complexity. This complexity is achieved
through two key ideas: the hierarchical low-rank structure of the matrix
is exploited to represent the original system as an extended sparse
system, and low-rank (re)compressions are applied to maintain the
sparsity pattern of the latter throughout the Gaussian elimination. The
solver has a tunable accuracy $\varepsilon$ allowing for a trade-off
between a highly accurate direct solver and a low-accurate
preconditioner.

Various numerical benchmarks have been carried out to validate the
accuracy and to assess the efficacy of the IFMM as a preconditioner. It
has been demonstrated that the IFMM preconditioner outperforms a
block-diagonal preconditioner for several types of boundary
conditions.  This is especially the case for pipe flow problems
involving many pores, which is relevant for various practical
applications.

Finally, it is important to note that the current IFMM implementation is
based on a uniform quadtree decomposition of the domain, which is far
from optimal for the complex geometries considered in this paper. The
use of an adaptive quadtree is expected to make the IFMM even more
efficient.

\section*{Acknowledgments}\label{sec:acks}
PC is a post-doctoral fellow of the Research Foundation Flanders (FWO)
and a Francqui Foundation fellow of the Belgian American Educational
Foundation (BAEF). The financial support is gratefully acknowledged.
BQ acknowledges startup funds from Florida State University.



\begin{thebibliography}{10}
\providecommand{\url}[1]{\texttt{#1}}
\providecommand{\urlprefix}{URL }
\expandafter\ifx\csname urlstyle\endcsname\relax
  \providecommand{\doi}[1]{doi:\discretionary{}{}{}#1}\else
  \providecommand{\doi}{doi:\discretionary{}{}{}\begingroup
  \urlstyle{rm}\Url}\fi

\bibitem{cue-jua2008}
Cueto-Felgueroso L, Juanes R. {Nonlocal Interface Dynamics and Pattern
  Formation in Gravity-Driven Unsaturated Flow Through Porous Media}.
  \emph{Physical Review Letters}  2008; \textbf{101}(24):244\,504.

\bibitem{kit-jia-val-chi-tsu-chr2014}
Kitamura K, Jiang F, Valocchi AJ, Chiyonobu S, Tsuji T, Christensen KT. {The
  study of heterogeneous two-phase flow around small-scale heterogeneity in
  porous sandstone by measured elastic wave velocities and lattice Boltzmann
  method simulation}. \emph{Journal of Geophysical Research: Solid Earth}
  2014; \textbf{119}(10):7564--7577.

\bibitem{mcd-hun-sit1986}
McDowell-Boyer LM, Hunt JR, Sitar N. {Particle Transport Through Porous Media}.
  \emph{Water Resources Research}  1986; \textbf{22}(13):1901--1921.

\bibitem{var-sta-pap1996}
Vardoulakis I, Stavropoulou M, Papanastasiou P. Hydromechanical aspects of the
  sand production problem. \emph{Transport in porous media}  1996;
  \textbf{22}(2):225--244.

\bibitem{wan2004}
Wang CY. {Fundamental Models for Fuel Cell Engineering}. \emph{Chemical
  Reviews}  2004; \textbf{104}(10):4727--4766.

\bibitem{lec-bot-wie2004}
Lecoanet HF, Bottero JY, Wiesner MR. {Laboratory Assessment of the Mobility of
  Nanomaterials in Porous Media}. \emph{Environmental Science \& Technology}
  2004; \textbf{38}(19):5164--5169.

\bibitem{wu-liu-xie-liu-sun2012}
Wu D, Liu H, Xie M, Liu H, Sun W. Experimental investigation on low velocity
  filtration combustion in porous packed bed using gaseous and liquid fuels.
  \emph{Experimental Thermal and Fluid Science}  2012; \textbf{36}:169--177.

\bibitem{wu-ye-sud2010}
Wu YS, Ye M, Sudicky EA. {Fracture-Flow-Enhanced Matrix Diffusion in Solute
  Transport Through Fractured Porous Media}. \emph{Transport in Porous Media}
  2010; \textbf{81}(1):21--34.

\bibitem{leb-den-vil2013}
Borgne TL, Dentz M, Villermaux E. {Stretching, Coalescence, and Mixing in
  Porous Media}. \emph{Physical Review Letters}  2013;
  \textbf{110}(20):204\,501.

\bibitem{hid-fe-cue-jua2012}
Hidalgo JJ, Fe J, Cueto-Felgueroso L, Juanes R. {Scaling of Convective Mixing
  in Porous Media}. \emph{Physical Review Letters}  2012;
  \textbf{109}(26):264\,503.

\bibitem{ber-sch-sil2000}
Berkowitz B, Scher H, Silliman SE. Anomalous transport in laboratory-scale,
  heterogeneous porous media. \emph{Water Resources Research}  2000;
  \textbf{36}(1):149--158.

\bibitem{sun-che-che2009}
Sun H, Chen W, Chen Y. Variable-order fractional differential operators in
  anomalous diffusion modeling. \emph{Physica A: Statistical Mechanics and its
  Applications}  2009; \textbf{388}(21):4586--4592.

\bibitem{hag-mck-mei2000}
Haggerty R, McKenna SA, Meigs LC. On the late-time behavior of tracer test
  breakthrough curves. \emph{Water Resources Research}  2000;
  \textbf{36}(12):3467--3479.

\bibitem{kaz-blo-kyr-chr2015}
Kazemifar F, Blois G, Kyritsis DC, Christensen KT. {Quantifying the flow
  dynamics of supercritical CO2--water displacement in a 2D porous micromodel
  using fluorescent microscopy and microscopic PIV}. \emph{Advances in Water
  Resources}  2015; .

\bibitem{szu-hes-jua2013}
Szulczewski ML, Hesse M, Juanes R. Carbon dioxide dissolution in structural and
  stratigraphic traps. \emph{Journal of Fluid Mechanics}  2013;
  \textbf{736}:287--315.

\bibitem{ica-boc-tem2016}
Icardi M, Boccardo G, Tempone R. {On the predictivity of pore-scale
  simulations: Estimating uncertainties with multilevel Monte Carlo}.
  \emph{Advances in Water Resources}  2016; .

\bibitem{mey-tch-jen2013}
Meyer DW, Tchelepi HA, Jenny P. A fast simulation method for uncertainty
  quantification of subsurface flow and transport. \emph{Water Resources
  Research}  2013; \textbf{49}(5):2359--2379.

\bibitem{den-leb-eng-bij2011}
Dentz M, Borgne TL, Englert A, Bijeljic B. {Mixing, spreading and reaction in
  heterogeneous media: A brief review}. \emph{Journal of Contaminant Hydrology}
   2011; \textbf{120}:1--17.

\bibitem{tel-dem-gre-kar2005}
T{\'e}l T, de~Moura A, Grebogi C, K{\'a}rolyi G. {Chemical and biological
  activity in open flows: A dynamical system approach}. \emph{Physics Reports}
  2005; \textbf{413}(2):91--196.

\bibitem{leb-oge-leb-hul-mar-sch1996}
Lebon L, Oger L, Leblond J, Hulin JP, Martys NS, Schwartz LM. {Pulsed gradient
  NMR measurements and numerical simulation of flow velocity distribution in
  sphere packings}. \emph{Physics of Fluids}  1996; \textbf{8}(2):293--301.

\bibitem{ica-boc-mar-tos-set2014}
Icardi M, Boccardo G, Marchisio DL, Tosco T, Sethi R. Pore-scale simulation of
  fluid flow and solute dispersion in three-dimensional porous media.
  \emph{Physical Review E}  2014; \textbf{90}(1):013\,032.

\bibitem{and-alm-fil-hav-suk-sta1997}
Jr JSA, Almeida MP, Filho JM, Havlin S, Suki B, Stanley HE. {Fluid Flow through
  Porous Media: The Role of Stagnant Zones}. \emph{Physical Review Letters}
  1997; \textbf{79}(20):3901.

\bibitem{blu-bij-don-gha-igl-mos-pal-pen2013}
Blunt MJ, Bijeljic B, Dong H, Gharbi O, Iglauer S, Mostaghimi P, Paluszny A,
  Pentland C. Pore-scale imaging and modelling. \emph{Advances in Water
  Resources}  2013; \textbf{51}:197--216.

\bibitem{man-gla-war1999}
Manz B, Gladden LF, Warren PB. {Flow and Dispersion in Porous Media:
  Lattice-Boltzmann and NMR Studies}. \emph{AIChE journal}  1999;
  \textbf{45}(9):1845--1854.

\bibitem{dea-leb-den-tar-bol-dav2013}
de~Anna P, Borgne TL, Dentz M, Tartakovsky AM, Bolster D, Davy P. {Flow
  Intermittency, Dispersion, and Correlated Continuous Time Random Walks in
  Porous Media}. \emph{Physical Review Letters}  2013;
  \textbf{110}(18):184\,502.

\bibitem{liu-che-lig1981}
Liu PL, Cheng AHD, Liggett JA, Lee JH. {Boundary Integral Equation Solutions to
  Moving Interface Between Two Fluids in Porous Media}. \emph{Water Resources
  Research}  1981; \textbf{17}(5):1445--1452.

\bibitem{lou-lee-kam1998}
Lough MF, Lee SH, Kamath J. {An Efficient Boundary Integral Formulation for
  Flow Through Fractured Porous Media}. \emph{Journal of Computational Physics}
   1998; \textbf{143}(2):462--483.

\bibitem{mal-gho-bir2014}
Malholtra D, Gholami A, Biros G. {A volume integral equation Stokes solver for
  problems with variable coefficients}. \emph{Proceedings of the International
  Conference for High Performance Computing, Networking, Storage and Analysis},
  IEEE Press, 2014; 92--102.

\bibitem{cam-ips-kel-mey-xue1996}
Campbell SL, Ipsen ICF, Kelley CT, Meyer CD, Xue ZQ. {Convergence estimates for
  solution of integral equations with GMRES}. \emph{Journal of Integral
  Equations and Applications}  1996; \textbf{8}(1):19--34.

\bibitem{gre-rok1987}
Greengard L, Rokhlin V. {A Fast Algorithm for Particle Simulations}.
  \emph{Journal of Computational Physics}  1987; \textbf{73}:325--348.

\bibitem{yin-bir-zor2004}
Ying L, Biros G, Zorin D. A kernel-independent adaptive fast multipole
  algorithm in two and three dimensions. \emph{Journal of Computational
  Physics}  2004; \textbf{194}(2):591--626.

\bibitem{dar-yor-ped1993}
Darden T, York D, Pedersen L. {Particle mesh Ewald: An $N \log (N)$ method for
  Ewald sums in large systems}. \emph{The Journal of Chemical Physics}  1993;
  \textbf{98}(12):10\,089--10\,092.

\bibitem{bar-hut1986}
Barnes J, Hut P. A hierarchical $\mathcal{O} ({N} \log {N})$ force-calculation
  algorithm. \emph{Nature Publishing Group}  1986; \textbf{324}:446--449.

\bibitem{qua-bir2015a}
Quaife B, Biros G. {On preconditioners for the Laplace double-layer potential
  in 2{D}}. \emph{Numerical Linear Algebra with Applications}  2015;
  \textbf{22}:101--122.

\bibitem{col-kre2012}
Colton D, Kress R. \emph{{Inverse Acoustic and Electromagnetic Scattering
  Theory}}, vol.~93. Springer Science \& Business Media, 2012.

\bibitem{hsi-wen2008}
Hsiao GC, Wendland WL. \emph{{Boundary Integral Equations}}. Springer: Berlin,
  2008.

\bibitem{heu:ste:tra1998}
Heuer N, Stephan EP, Tran T. {Multilevel additive Schwarz method for the h-p
  version of the Galerkin boundary element method}. \emph{Mathematics of
  Computation of the American Mathematical Society}  1998;
  \textbf{67}(222):501--518.

\bibitem{mai-ste-tra2000}
Maischak M, Stephan EP, Tran T. {Multiplicative Schwarz algorithms for the
  Galerkin boundary element method}. \emph{SIAM Journal on Numerical Analysis}
  2000; \textbf{38}(4):1243--1268.

\bibitem{fun-ste1997}
Funken SA, Stephan EP. {The BPX preconditioner for the single layer potential
  operator}. \emph{Applicable Analysis}  1997; \textbf{67}(3-4):327--340.

\bibitem{che2000}
Chen K. {An Analysis of Sparse Approximate Inverse Preconditioners for Boundary
  Integral Equations}. \emph{SIAM Journal on Matrix Analysis and Applications}
  2000; \textbf{22}(4):1058--1078.

\bibitem{gra-kum-sam1998}
Grama A, Kumar V, Sameh A. {Parallel Hierarchical Solvers and Preconditioners
  for Boundary Element Methods}. \emph{SIAM Journal on Scientific Computing}
  1997; \textbf{20}(1):337--358.

\bibitem{beb2005}
Bebendorf M. {Hierarchical LU decomposition-based preconditioners for BEM}.
  \emph{Computing}  2005; \textbf{74}(3):225--247.

\bibitem{nab-kor-lei-whi1994}
Nabors K, Korsmeyer FT, Leighton FT, White JK. Preconditioned, adaptive,
  multipole-accelerated iterative methods for three-dimensional first-kind
  integral equations of potential theory. \emph{SIAM Journal on Scientific and
  Statistical Computing}  1994; \textbf{15}:713--735.

\bibitem{vav1992}
Vavasis SA. Preconditioning for boundary integral equations. \emph{SIAM Journal
  on Matrix Analysis and Applications}  1992; \textbf{13}(3):905--925.

\bibitem{ste-wen1998}
Steinbach O, Wendland W. The construction of some efficient preconditioners in
  the boundary element method. \emph{Advances in Computational Mathematics}
  1998; \textbf{9}:191--216.

\bibitem{bra-ley-pas1994}
Bramble JH, Leyk Z, Pasciak JE. {The Analysis of Multigrid Algorithms for
  Pseudodifferential Operators of Order Minus One}. \emph{Mathematics of
  Computation}  October 208; \textbf{63}(208):461--478.

\bibitem{hsi-xu-zha2014}
Hsiao GC, Xu L, Zhang S. {Solving Negative Order Equations by the Multigrid
  Method Via Variable Substitution}. \emph{Journal of Scientific Computing}
  2014; \textbf{59}(2):371--385.

\bibitem{lan-pus-rei2003}
Langer U, Pusch D, Reitzinger S. Efficient preconditioners for boundary element
  matrices based on grey-box algebraic multigrid methods. \emph{International
  Journal for Numerical Methods in Engineering}  2003;
  \textbf{58}(13):1937--1953.

\bibitem{of2008}
Of G. An efficient algebraic multigrid preconditioner for a fast multipole
  boundary element method. \emph{Computing}  2008; \textbf{82}:139--155.

\bibitem{lan-pus2005}
Langer U, Pusch D. Data-sparse algebraic multigrid methods for large scale
  boundary element equations. \emph{Applied Numerical Mathematics}  2005;
  \textbf{54}(3):406--424.

\bibitem{amin14a}
Aminfar A, Ambikasaran S, Darve E. {A Fast Block Low-Rank Dense Solver with
  Applications to Finite-Element Matrices}. \emph{arXiv}  2014; :1403.5337.

\bibitem{chan06a}
Chandrasekaran M S~Gu, Pals T. A fast {U}{L}{V} decomposition solver for
  hierarchically semi-separable representations. \emph{SIAM Journal on Matrix
  Analysis and Applications}  2006; \textbf{28}(3):603--622.

\bibitem{shen07b}
Sheng Z, Dewilde P, Chandrasekaran S. {Algorithms to Solve Hierarchically
  Semi-separable Systems}. \emph{System Theory, the Schur Algorithm and
  Multidimensional Analysis}, \emph{Operator Theory: Advances and
  Applications}, vol. 176, Alpay D, Vinnikov V (eds.). Birkh\"{a}user Basel,
  2007; 255--294.

\bibitem{bebe08a}
Bebendorf M. \emph{{Hierarchical Matrices: A Means to Efficiently Solve
  Elliptic Boundary Value Problems}}. 1st edn., Springer Publishing Company,
  2008.

\bibitem{borm2007data}
B{\"o}rm S. Data-sparse approximation of non-local operators by
  $\mathcal{H}^2$-matrices. \emph{Linear Algebra and its Applications}  2007;
  \textbf{422}(2):380--403.

\bibitem{hackbusch2002data}
Hackbusch W, B{\"o}rm S. {Data-sparse Approximation by Adaptive
  $\mathcal{H}^2$-Matrices}. \emph{Computing}  2002; \textbf{69}(1):1--35.

\bibitem{gill14a}
Gillman A, Martinsson PG. An ${O}({N})$ algorithm for constructing the solution
  operator to 2{D} elliptic boundary value problems in the absence of body
  loads. \emph{Advances in Computational Mathematics}  2014;
  \textbf{40}(4):773--796.

\bibitem{kong11a}
Kong W, Bremer J, Rokhlin V. An adaptive fast direct solver for boundary
  integral equations in two dimensions. \emph{Applied and Computational
  Harmonic Analysis}  2011; \textbf{31}(3):346--369.

\bibitem{mart05b}
Martinsson PG, Rokhlin V. A fast direct solver for boundary integral equations
  in two dimensions. \emph{Journal of Computational Physics}  2005;
  \textbf{205}(1):1--23.

\bibitem{mar-bar-gil-vee2015}
Marple GR, Barnett A, Gillman A, Veerapaneni S. A fast algorithm for simulating
  multiphase flows through periodic geometries of arbitrary shape. \emph{arXiv}
   2015; \textbf{1510.05616}.

\bibitem{gil-bar2013}
Gillman A, Barnett A. A fast direct solver for quasi-periodic scattering
  problems. \emph{Journal of Computational Physics}  2013;
  \textbf{248}:309--322.

\bibitem{gil-bar-mar2014}
Gillman A, Barnett A, Martinsson PG. A spectrally accurate direction solution
  technique for frequency-domain scattering problems with variable media.
  \emph{BIT Numerical Mathematics}  2014; \textbf{55}(1):141--170.

\bibitem{amb-dar2014}
Ambikasaran S, Darve E. {The Inverse Fast Multipole Method}. \emph{arxiv}
  2014; \textbf{1407.1572}.

\bibitem{ij-sjsc-coul-a}
Coulier P, Pouransari H, Darve E. The inverse fast multipole method: using a
  fast approximate direct solver as a preconditioner for dense linear systems.
  \emph{arXiv}  2015; :1508.01\,835.

\bibitem{pow1993}
Power H. {The completed double layer boundary integral equation method for
  two-dimensional Stokes flow}. \emph{IMA Journal of Applied Mathematics}
  1993; \textbf{51}(2):123--145.

\bibitem{poz1992}
Pozrikidis C. \emph{{Boundary Integral and Singularity Methods for Linearized
  Viscous Flow}}. Cambridge University Press: New York, NY, USA, 1992.

\bibitem{kapu97a}
Kapur S, Rokhlin V. {High-Order Corrected Trapezoidal Quadrature Rules for
  Singular Functions}. \emph{SIAM Journal on Numerical Analysis}  1997;
  \textbf{34}(4):1331--1356.

\bibitem{alp1999}
Alpert BK. {Hybrid Gauss-Trapezoidal Quadrature Rules}. \emph{SIAM Journal on
  Scientific Computing}  1999; \textbf{20}:1551--1584.

\bibitem{bar-wu-vee2014}
Barnett A, Wu B, Veerapaneni S. {Spectrally-Accurate Quadratures for Evaluation
  of Layer Potentials Close to the Boundary for the 2{D} Stokes and Laplace
  Equations}. \emph{arXiv}  2014; \textbf{1410.2187}.

\bibitem{yin-bir-zor2006}
Ying L, Biros G, Zorin D. A high-order 3{D} boundary integral equation solver
  for elliptic {PDE}s in smooth domains. \emph{Journal of Computational
  Physics}  2006; \textbf{219}(1):247--275.

\bibitem{cou-dar2016}
Coulier P, Darve E. {Efficient mesh deformation based on radial basis function
  interpolation by means of the inverse fast multipole method}. \emph{Computer
  Methods in Applied Mechanics and Engineering}  2016; .

\bibitem{xia10a}
Xia J, Chandrasekaran S, Gu M, Li X. Fast algorithms for hierarchically
  semiseparable matrices. \emph{Numerical Linear Algebra with Applications}
  2010; \textbf{17}(6):953--976.

\bibitem{xia10b}
Xia J, Chandrasekaran S, Gu M, Li X. {Superfast Multifrontal Method for Large
  Structured Linear Systems of Equations}. \emph{SIAM Journal on Matrix
  Analysis and Applications}  2010; \textbf{31}(3):1382--1411.

\bibitem{ambi13b}
Ambikasaran S, Darve E. {An $\mathcal{O} ({N} \log {N})$ Fast Direct Solver for
  Partial Hierarchically Semi-Separable Matrices}. \emph{Journal of Scientific
  Computing}  2013; \textbf{57}(3):477--501.

\bibitem{coro15a}
Corona E, Martinsson PG, Zorin D. A ${O}({N})$ direct solver for integral
  equations on the plane. \emph{Applied and Computational Harmonic Analysis}
  2015; \textbf{38}(2):284--317.

\bibitem{lize14a}
Liz\'{e} B. {Fast Direct Solver for the Boundary Element Method in
  Electromagnetism and Acoustics: $\mathcal{H}$-Matrices. Parallelism and
  Industrial Applications}. Ph{D} {T}hesis, Universit\'{e} Paris 13 2014.

\bibitem{saa-sch1986}
Saad Y, Schultz M. {GMRES: A Generalized Minimal Residual Algorithm for Solving
  Nonsymmetric Linear Systems}. \emph{SIAM Journal on Scientific and
  Statistical Computing}  1986; \textbf{7}:856--869.

\bibitem{fong09a}
Fong W, Darve E. The black-box fast multipole method. \emph{Journal of
  Computational Physics}  2009; \textbf{228}(23):8712--8725.

\bibitem{hsi-kop-wen1980}
Hsiao GC, Kopp P, Wendland WL. {A Galerkin Collocation Method for Some Integral
  Equations of the First Kind}. \emph{Computing}  1980; \textbf{25}(2):89--130.

\bibitem{hal2001}
Haller G. Distinguished material surfaces and coherent structures in
  three-dimensional fluid flows. \emph{Physica D}  2001; \textbf{149}:248--277.

\end{thebibliography}

\end{document}